\documentclass[10pt,a4paper]{article}
\usepackage[utf8]{inputenc}
\usepackage[T1]{fontenc}
\usepackage{amssymb}
\usepackage{amsthm}
\usepackage{color}
\definecolor{DarkOlive}{rgb}{0.1047,0.2412,0.0064}
\definecolor{FireBrick}{rgb}{0.5812,0.0074,0.0083}
\definecolor{RoyalBlue}{rgb}{0.0236,0.0894,0.6179}
\definecolor{RoyalGreen}{rgb}{0.0236,0.6179,0.0894}
\definecolor{RoyalRed}{rgb}{0.6179,0.0236,0.0894}
\definecolor{LightBlue}{rgb}{0.8544,0.9511,1.0000}
\definecolor{Black}{rgb}{0.0,0.0,0.0}
\definecolor{MidnightBlue}{rgb}{0.0035,0.0020,0.1363}
\definecolor{FireBrick3}{rgb}{0.5812,0.0074,0.0083}
\definecolor{FireBrick4}{rgb}{0.2156,0.0023,0.0035}
\definecolor{Blue2}{rgb}{0.0000,0.0000,0.8644}
\definecolor{Navy}{rgb}{0.0000,0.0000,0.1927}
\definecolor{MediumBlue}{rgb}{0.0000,0.0000,0.6179}
\usepackage[
        bookmarks=false,pdftitle={Greatest common divisors of values of 
        integer polynomials and an application to maximal tori},
        colorlinks=true,backref=false,breaklinks=true,linkcolor=MediumBlue,
        citecolor=FireBrick3,filecolor=RoyalRed,
        urlcolor=Blue2]{hyperref}

\theoremstyle{plain}
\newtheorem{Thm}{Theorem}[section]

\newtheorem{Pro}[Thm]{Proposition}

\theoremstyle{definition}

\newtheorem{Rem}[Thm]{Remark}
\newtheorem{Alg}[Thm]{Algorithm}
\newtheorem{Exe}[Thm]{Example}

\newenvironment{Prf}{\noindent\textbf{Proof.}}{\hfill $\Box$ \medskip}

\newcommand{\mod}{ \textrm{ mod } }

\newcommand{\F}{\mathbb{F}}
\newcommand{\Z}{\mathbb{Z}}

\newcommand{\Q}{\mathbb{Q}}
\newcommand{\GL}{\textrm{GL}}
\newcommand{\bG}{\ensuremath{G}}
\newcommand{\GAP}{{GAP}}
\newcommand{\CHEVIE}{{CHEVIE}}

\newcommand{\lcm}{{\rm lcm}}

\begin{document}
% MSC 11A05 20D06 20G40 (Euclidean algorithm, simple groups, groups of Lie
% type
\title{Greatest common divisors of values of integer polynomials and
an application to maximal tori}
\author{Frank Lübeck\thanks{This is a contribution to 
Project-ID 286237555 -- TRR 195 -- by the
Deutsche Forschungsgemeinschaft (DFG, German Research Foundation)}}
\maketitle
\begin{abstract}
We describe how to compute for two polynomials $f(X), g(X) \in \Z[X]$ with
integer coefficients the greatest common divisors of $f(z)$ and $g(z)$ for
all integers $z \in \Z$.

As an application we determine the structures (elementary divisors) of all
maximal tori in exceptional groups of Lie type.
\end{abstract}

\section{Greatest commom divisors}\label{secgcdzx}
Given polynomials $f(X), g(X) \in \Z[X]$ over the integers, we describe an
algorithm to  determine the greatest common  divisor $\gcd(f(z),g(z))$ for
the evaluations of $f$ and $g$ at all integers $z \in \Z$.

The algorithm  will be a  constructive proof of the  following proposition
which also suggests how to present the result.

\begin{Pro} \label{gengcd}
Let $f(X), g(X) \in \Z[X]$. 
Then there exists an integer $m \in \Z_{>0}$ such that:
\begin{itemize}
\item[(i)]
For each $i \in \Z$, $0 \leq i < m$, there exists a polynomial $h_i(X) \in
\Z[X]$ such  that $\gcd(f(z),g(z)) =  h_i(z)$ for all  $z \in \Z$  with $z
\equiv  i \mod  m$. The  polynomials $h_i(X)$  for all  $i$ have  the same
primitive part.
\item[(ii)]
For each $i \in \Z$, $0 \leq i < m$, there exists two polynomials $a_i(X),
b_i(X)  \in \Q[X]$  such that  $a_i(X) f(X)  + b_i(X)  g(X) =  h_i(X)$ and
$a_i(z), b_i(z) \in \Z$ for all $z \in \Z$ with $z \equiv i \mod m$.
\end{itemize}
\end{Pro}

Recall that any polynomial $0 \neq f(X) \in \Q[X]$ can be uniquely written
as $f(X)  = c(f)  \cdot \tilde{f}(X)$  where $c(f) \in  \Q$ is  called the
\emph{content} of  $f$ and  $\tilde{f}(X) \in \Z[X]$  is \emph{primitive},
that is  the greatest common  divisor of its  coefficients is $1$  and its
leading coefficient is positive.

There are  polynomials $f(X) \in  \Q[X]$ such that  $f(z) \in \Z$  for all
integers  or all  integers  in certain  congruence  classes. The  greatest
common divisors of such values can also be investigated with the algorithm
described here by first  substiting $X = kY+i$ for $0 \leq  i <k$ were $k$
is a common denominator of the coefficients of $f$.

Here is the algorithm, we include some explanations.

\begin{Alg} \label{alggcd}\mbox{}

\textbf{Input:} $f(X), g(X) \in \Z[X]$.

\textbf{Output:} $m$ and the $h_i$, $a_i$ and $b_i$ as in
Proposition~\ref{gengcd}.

\textbf{Step 1.} [primitive part] 

\noindent
Compute  the greatest  common divisor  of $f$  and $g$  as polynomials  in
$\Q[X]$ with  the Euclidean  algorithm. Its primitive  part $h(X)$  is the
primitive part of all $h_i(X)$.

\noindent
Substitute $f(X)$  by $f(X)/h(X)$  and $g(X)$ by  $g(X)/h(X)$. So,  in the
following  steps  $f$ and  $g$  are  relatively  prime as  polynomials  in
$\Q[X]$.

\noindent
[\emph{Note  that we  still  have  $f(X), g(X)  \in  \Z[X]$. This  follows
from  Gauss's  lemma  which  shows  that  the  content  of  a  product  of
polynomials  in  $\Q[X]$  is  the  product  of  their  contents,  see  for
example~\cite[2.16]{JaI}.}]

\textbf{Step 2.} [find $m$ as linear combination]

\noindent
Use the extended Euclidean algorithm to compute polynomials $\tilde{a}(X),
\tilde{b}(X) \in \Q[X]$ with $\tilde{a}(X)  f(X) + \tilde{b}(X) g(X) = 1$.
Let $m \in \Z_{>0}$ be the least common denominator of all coefficients of
$\tilde{a}$  and $\tilde{b}$.  Then set  $a(X) :=  m \cdot  \tilde{a}(X)$,
$b(X) := m  \cdot \tilde{b}(X) \in \Z[X]$, such that  we get the following
equation in $\Z[X]$:
\[ {a}(X) f(X) + {b}(X) g(X) = m. \]

\noindent
[\emph{This shows that  for any $z \in \Z$ the  greatest common divisor of
$f(z)$ and $g(z)$ is  a divisor of $m$. Since $f(z+m)  \equiv f(z) \mod m$
and  $g(z+m)  \equiv  g(z)  \mod  m$ this  divisor  only  depends  on  the
congruence class of $z \mod m$.}]

\textbf{Step 3.} [find the $h_i$, $a_i$, $b_i$]

\noindent
For each $i \in \Z$ with $0 \leq i < m$ do:

\quad Set $c_i := \gcd(f(i), g(i))$ and $h_i := c_i \cdot h$ (the $h$ from
step 1).

\quad Set $r_i := m/c_i \in \Z$ and  find $x_i \in \Z$, $0 \leq x_i < r_i$
such that

\quad\quad $a(i) - x_i \frac{g(i)}{c_i} \equiv 0 \mod r_i$ and $b(i) + x_i
\frac{f(i)}{c_i} \equiv 0 \mod r_i$.

\quad  Set  $a_i(X)  :=  \frac{1}{r_i}(a(X)  -\frac{x_i}{c_i}  g(X))$  and
$b_i(X) := \frac{1}{r_i}(b(X) + \frac{x_i}{c_i} f(X))$.

\noindent
[\emph{Note  that  $f(z)/c_i$ and  $g(z)/c_i$  are  integers in  the  same
congruence class modulo $r_i$ for all $z \in \Z$ with $z \equiv i \mod m$.
Also, $a(z)  \equiv a(i)  \mod m$ and  $b(z) \equiv b(i)  \mod m$  (and so
modulo $r_i$) for all $z \equiv i  \mod m$. So, the determination of $x_i$
above would yield  the same result if  we used the evaluation  at any such
$z$ instead of  the evaluation at $z  = i$.\\ It is  clear by construction
that for all $z  \in \Z$ with $z \equiv i \mod m$  we have $a_i(z), b_i(z)
\in \Z$ and
\[ a_i(z) f(z) + b_i(z) g(z) = \frac{m}{r_i} = c_i.\]}]

\textbf{Step 4.} [find smaller $m$, return result]

\noindent
Find the  smallest divisor  $d$ of  $m$ from step~2  such that  the $c_i$,
$a_i(X)$, $b_i(X)$  in step~3 only  depend on  the congruence class  of $i
\mod d$. Then redefine $m := d$.

Return $m$ and  the $h_i(X), a_i(X), b_i(X)$  from step~3 for $0  \leq i <
m$.

\noindent
[\emph{It  may be  more efficient  to return  $h(X)$ from  step~1 and  the
contents $c_i$ from step~3 instead of the $h_i(X)$.}]
\end{Alg}
\begin{Prf}
We will show that in step~3 there  exists an $x_i \in \Z$ with the desired
property. That the algorithm terminates and  is correct is then clear from
the comments given within the algorithm.

So, let $F,G,A,A',B,B',M \in \Z$ with  $AF+BG = A'F+B'G = M$, $\gcd(F,G) =
C$ (so $C \mid M$) and $R = \frac{M}{C}$. Then we have $(A-A')F = (B'-B)G$
and since  this number is a  multiple of $F$ and  $G$ it is a  multiple of
$\lcm(F,G)  = \frac{F  G}{C}$, hence  $(A-A') =  x \cdot  \frac{G}{C}$ and
$(B-B') = x \cdot \frac{F}{C}$ for some  $x \in \Z$. This shows how to get
from one  linear combination $AF  + BG = M$  all other pairs  $A',B'$ with
$A'F + B'G =  M$. It remains to show that there exist  such $A', B'$ which
are  both  divisible by  $R$.  But  this  is  clear because  the  extended
Euclidean algorithm yields  $I,J \in \Z$ with  $IF + JG = C$,  so $(IR)F +
(JR)G = CR = M$. Since $A - x G /  C \equiv A - (x+R) G / C \mod R$ we can
find an appropriate $x$ in the range $0 \leq x < R$.
\end{Prf}

\begin{Exe}\mbox{} \label{gcdex}

\begin{itemize}
\item[(a)] $f(X) = X - 1$ and $g(X) = X + 1$. Here we have $-f(X) + g(X) =
2$, so that the  greatest common divisor of $f(z)$ and  $g(z)$ will be $1$
or $2$ for all integers $z  \in \Z$. Algorithm~\ref{alggcd} yields in this
case:

$m = 2$.

If $z \equiv  0 \mod 2$ then $\gcd(f(z),  g(z)) = 1$ and $a_0(X)  = -X/2 -
1$, $b_0(X) = X/2$. (So, for even  $z$ we have $a_0(z), b_0(z) \in \Z$ and
$a_0(z) f(z) + b_0(z) g(z) = 1$.)

If $z  \equiv 1 \mod  2$ then  $\gcd(f(z), g(z)) =  2$, $a_1(X) =  -1$ and
$b_1(X) = 1$.
\item[(b)]
$f(X) = X^2+X+1$ and  $g(X) = 2 X^2+2$. In this  case, our algorithm shows
that  $\gcd(f(z), g(z))  = 1$  for all  $z  \in \Z$.  So, if  we are  only
interested in the greatest  common divisor, we have $m =  1$ and $h_0(X) =
1$. But it  is not difficult to  see that there are  no polynomials $a(X),
b(X) \in \Q[X]$  such that $a(X)f(X) +  b(X)g(X) = 1$ and  $a(z), b(z) \in
\Z$ for all $z \in \Z$. Our algorithm~\ref{alggcd} yields:

$m = 2$.

$z  \equiv  0 \mod  2$:  $h_0(X)  = 1$,  $a_0(X)  =  -X^2-X-1$, $b_0(X)  =
X^2/2+X+1$.

$z \equiv 1 \mod 2$: $h_1(X) = 1$, $a_1(X) = -X$, $b_1(X) = (X+1)/2$.

\item[(c)] 
$f(X) := 5X^{4}+12X^{3}+2X^{2}+3X+14$ and $g(X) := 3X^{3}-10X+4$.

This example  shows that  the module  $m$ in  Proposition~\ref{gengcd} can
easily  become  large,  even  for  polynomials  which  do  not  look  very
complicated. In this case our algorithm finds $m = 5291$, and depending on
the congruence class modulo $m$ we have
\[\gcd(f(z), g(z)) \in \{1, 11, 13, 37, 143, 407, 481, 5291\}.\] 
Here the polynomials $a_i(X), b_i(X)$ are pairwise different for all cases
$0 \leq i < 5291$.
\end{itemize}
\end{Exe}

\begin{Rem}
Depending on  implementation details of  the used Euclidean  algorithm the
module $m$ found in step~2  of Algorithm~\ref{alggcd} may be unnecessarily
large. Sometimes we were successful finding  a smaller $m$ by applying the
Euclidean algorithm  to several unimodular integer  linear combinations of
$f$ and  $g$. If one  finds different $m$ one  can combine the  results to
find a linear  combination of $f$ and $g$ for  the greatest common divisor
of all those $m$.

\end{Rem}

\section{Application: Structure of maximal tori in exceptional groups}

\subsection{Setup}

In this  section we consider  certain finite groups  of Lie type.  We only
give a very  brief survey of the  setup, see some standard  texts for more
details (e.g.,  \cite[Ch. 1  and 3]{Ca85}  for the  theoretical background
and~\cite[Sec.2]{BL13} for the computational setup).

These groups  can be  described in  terms of  a reductive  algebraic group
$\bG$ over the algebraic closure $\bar{\F}_p$  of a finite prime field and
a Frobenius morphism  $F$ of $\bG$ associated  to a power $q$  of $p$. The
finite group of fixed  points $G(q) := \bG^F$ is called  a finite group of
Lie type.

An important  conjugacy class  of subgroups of  $\bG$ consists  of maximal
tori $T$. Similarly,  maximal tori of $\bG^F$ are subgroups  of type $T^F$
for maximal tori $T$ with $F(T) = T$.

The algebraic  group $\bG$  is determined  up to  isomorphism by  its root
datum with respect to any maximal torus $T \subset \bG$ and the underlying
field. The root datum consists of  a $\Z$-lattice $X \cong \Z^r$, its dual
lattice  $Y \cong  \Z^r$  and  a bijection  between  finite subsets  $\Phi
\subset X$ and $\Phi^\vee \subset Y$, called the roots and coroots.

$\Phi$  and $\Phi^\vee$  can be  constructed from  a linearly  independent
subset  of   simple  roots   $\Delta  =  \{\alpha_1,   \ldots,  \alpha_l\}
\subset  \Phi$ and  its  coroots $\Delta^\vee  = \{\alpha_1^\vee,  \ldots,
\alpha_l^\vee\} \subset  \Phi^\vee$. For each $\alpha_i  \in \Delta$ there
is  a reflection  on  $Y$  defined by  $s_i:  Y \to  Y$,  $y  \mapsto y  -
y(\alpha_i)  \alpha_i^\vee$.  The elements  $S  :=  \{s_1, \ldots,  s_l\}$
generate the Weyl group  $W$ of the root datum and $\bG$ and  $S$ is a set
of Coxeter  generators of $W$.  The matrix $C  = (c_{ij})$ with  $c_{ij} =
\alpha_i^\vee(\alpha_j)$ is  called the Cartan  matrix of the  root system
and it can be encoded by a  Dynkin diagram. When the maximal torus $T$ for
the root datum is $F$-invariant, then  the Frobenius $F$ induces a natural
$\Z$-linear map on $Y$ which has the form $q F_0$ where $F_0$ is of finite
order. Also, this induces an automorphism  on $W$. One can choose $T$ such
that $\Delta^\vee$ and $S$ are invariant under $F$.

The $\bG^F$-conjugacy classes of maximal tori of $\bG^F$ are parameterized
by the $F$-conjugacy classes of the Weyl group $W$.

Let $\Q_{p'}$  be the additive  subgroup of  the rational numbers  $\Q$ of
numbers with denominator not divisible by $p$. Then we have
\[ T \cong Y \otimes (\Q_{p'}/\Z) \]
and when  $w \in W$  then the action of  the Frobenius on  a corresponding
maximal  torus is  equivalent to  the action  of $w  F$ on  $Y$, naturally
extended to $T = Y \otimes (\Q_{p'}/\Z)$.

\subsection{Simply-connected exceptional groups}

In this paper we are mainly interested in the case where $\bG$ is a simple
group of exceptional  type. We first consider  simply-connected groups. In
this case we can choose (by definition) the simple coroots as a $\Z$-basis
of $Y$. The simple  roots written in the dual basis are  then given by the
rows of  the transposed  Cartan matrix  of the root  system of  $\bG$. The
finite  order map  $F_0$  is then  given by  the  permutation matrix  that
permutes the simple  coroots according to the automorphism  induced by $F$
on the  Dynkin diagram (so just  the identity for untwisted  groups). With
respect to the chosen basis we  identify $T = Y \otimes (\Q_{p'}/\Z)$ with
$(\Q_{p'}/\Z)^r$.

For a maximal torus $T_w$ in the $\bG^F$-conjugacy class associated to the
$F$-conjugacy class of $w \in W$ we  get that $T_w^F$ is isomorphic to the
group of the $t \in (\Q_{p'}/\Z)^r$ with
\[t (q w_Y F_0 - 1) = 0\] 
(here $w_Y$ is the matrix of $w$  acting on $Y$ with respect to the chosen
basis, and $1$ is the identity matrix).

For a  fixed prime  power $q$ we  have $M(q) :=  q w_Y  F_0 - 1  \in \Z^{r
\times r}$  and we can  find the solutions  as follows: Compute  the Smith
normal form of $M(q)$ with transforming matrices $L, R \in \GL_r(\Z)$ such
that $D =  L M(q) R$ is  a diagonal matrix with entries  $d_1, \ldots, d_r
\in \Z$  (the elementary divisors  of $M(q)$). Then $t  M(q) = 0$  has the
same set of solutions as $t M(q) R = 0$ which we can write as
\[ t M(q) R = t L^{-1} (L  M(q) R) = (t L^{-1}) D = 0. \]
Note that $q$ is a power of $p$ and $M(q) \mod p$ is $-1$, so $\det(M(q))$
and the $d_1, \ldots, d_r$ are not divisible by $p$. From this we see that
\[t L^{-1} = (i_1/d_1, \ldots, i_r/d_r) \mod Y \in
(\Q_{p'}/\Z)^r \textrm{ with } i_1, \ldots i_r \in \Z.\] 
So, the  structure of $T_w^F$  is the direct  product of cyclic  groups of
orders given by the elementary divisors of $M(q)$.

Instead of computing  with a matrix $M(q)$ for a  specific prime power $q$
we can use the results from Section~\ref{secgcdzx} and compute with $M(q)$
where $q$ is considered as indeterminate  over $\Z$, so $M(q) \in \Z[q]^{r
\times  r}$. We  can  compute  the product  of  the  first $i$  elementary
divisors,  $d_1\cdots  d_i$,  by  computing the  greatest  common  divisor
of  all $(i  \times  i)$-minors of  $M(q)$, using  algorithm~\ref{alggcd}.
From proposition~\ref{gengcd}  we know  that the  result can  be described
uniformly  for all  $q$ in  a fixed  congruence class  of $q$  modulo some
integer $m$.

For the  computations we used an  implementation of algorithm~\ref{alggcd}
in the system \GAP~\cite{GAP}. The  matrices $M(q)$ are straightforward to
compute with \CHEVIE~\cite{CHEVIE}.

\subsection{Remarks}\label{secrem}
\begin{enumerate}
\item \textbf{Parameterization.}  Instead of just computing  the structure
of $T_w^F$  we may also want  to find an explicit  parameterization of the
solutions of $t M(q) = 0$. In  principle this can be done by computing the
Smith normal form  of $M(q)$ with transforming matrices  $L(q)$ and $R(q)$
with  $q$ as  a  parameter, using~\ref{gengcd}(ii).  But  in practice  our
implementation of such an algorithm did  not work well because it may lead
to a  large modulus $m$  and too many cases,  while a computation  by some
easy heuristics or by hand sometimes finds solutions which are independent
of congruence conditions on $q$.
\item \textbf{More uniform description.} It  can happen that the structure
of $T_w^F$  as some product of  cyclic groups can be  described uniformely
for  all $q$  although the  elementary divisors  may depend  on congruence
conditions on  $q$. A typical  example would be a  torus that is  a direct
product of a cyclic group of order $q-1$ and a cyclic group of order $q+1$
(for all  $q$). But this  leads to two  cases for the  elementary divisors
(see example~\ref{gcdex}(a)): If $q$ is even  then this group is cyclic of
order  $q^2-1$. And  if $q$  is  odd then  the  group is  not cyclic,  the
elementary divisors are $2$ and $(q^2-1)/2$.
\item \textbf{Adjoint groups.}  If $\bG$ is a simple  exceptional group of
adjoint type  then its  dual group is  isomorphic to  the simply-connected
group of the  same type. A torus  $T_w^F$ in $\bG^F$ is  isomorphic to the
character group of  the corresponding torus in the dual  group, and so the
structure as  abelian group  is the  same. See~\cite[Ch.4,  4.4]{Ca85} for
more details.
\item   \textbf{Finite   simple  groups.}   In   many   cases  of   simple
simply-connected groups $\bG$  of exceptional type the group  $\bG^F$ is a
finite simple group.  The exceptions are type $E_6$ with  $q \equiv 1 \mod
3$, ${}^2E_6$  with $q \equiv 2  \mod 3$ and  type $E_7$ with $q  \equiv 1
\mod 2$. In these  cases the center $Z$ of $G(q)$ is of  order 3, 3 and 2,
respectively, $Z$  is contained in  any maximal  torus, and $G(q)/Z$  is a
finite simple group.

We do also compute the structure of $T_w^F/Z$ in these cases. 

The center $Z$ consists of the solutions of  $t C = 0$ in $T$ where $C$ is
the  Cartan matrix  of the  root system  ($Z$ is  the intersection  of the
kernels  of all  roots). We  first find  transforming matrices  $L, R  \in
\GL_r(\Z)$ such that  $L C R$ has  Smith normal form. Then  the entries in
the last column of $L M(q)$ are divisible by the last diagonal entry of $L
C R$ (that is $2$ or $3$,  respectively). Dividing them yields a system of
equations for $T_w^F/Z$.
\item \textbf{Suzuki and Ree groups.} With  a small adjustment we can also
handle the Suzuki  and Ree groups with our methods.  These occur for $\bG$
of type $B_2$,  $G_2$ or $F_4$ and an exceptional  Frobenius morphism that
exists  only in  characteristic $p  = 2$,  $3$ or  $2$, respectively.  The
parameter $q$ is then an odd power of $\sqrt{p}$.

Here we substitute $F_0$ by  a monomial matrix with entries $\sqrt{p}^{\pm
1}$ and  substitute $q$  by the  integer $q' =  q/\sqrt{p}$. This  way the
equations for  $T_w^F$ involve  matrices in $\Z[q']$  which we  can handle
with the methods in Section~\ref{secgcdzx}.
\end{enumerate}

\section{Tables}

For each type of group we give a table with the following information:
\begin{itemize}
\item the type,
\item  names   of  $F$-conjugacy  classes   of  $W$,  as  in   the  tables
of~\cite{CaWeyl72},
\item reduced words of representatives  $w \in W$ of $F$-conjugacy classes
(we write $14321$ instead of $s_1 s_4 s_3 s_2 s_1$),
\item the structure of the torus $T_w^F$ as list of non-trivial elementary
divisors (we write $\Phi_i$ for the $i$-th cyclotomic polynomial evaluated
at $q$, e.g. $\Phi_1 = q-1$, $\Phi_8 = q^4+1$, see~\cite[13.9]{Ca85} for a
list),
\item in cases where the structure  depends on the congruence class of $q$
we also give a uniform description of  the structure of $T_w^F$ as list of
orders of cyclic factors (see~\ref{secrem} 2.).
\end{itemize}

The simple roots and  generators of the Weyl group are  numbered as in the
following Dynkin diagrams:

\begin{center}
{
\newlength{\tmplen}
\setlength{\tmplen}{\unitlength}
\setlength{\unitlength}{0.8pt}
\scriptsize
\begin{picture}(360,190)
\put( 10, 40){$E_7$}
\put( 40, 40){\circle*{5}}
\put( 38, 45){1}
\put( 40, 40){\line(1,0){20}}
\put( 60, 40){\circle*{5}}
\put( 58, 45){3}
\put( 60, 40){\line(1,0){20}}
\put( 80, 40){\circle*{5}}
\put( 78, 45){4}
\put( 80, 40){\line(0,-1){20}}
\put( 80, 20){\circle*{5}}
\put( 85, 18){2}
\put( 80, 40){\line(1,0){20}}
\put(100, 40){\circle*{5}}
\put( 98, 45){5}
\put(100, 40){\line(1,0){20}}
\put(120, 40){\circle*{5}}
\put(118, 45){6}
\put(120, 40){\line(1,0){20}}
\put(140, 40){\circle*{5}}
\put(138, 45){7}

\put(190, 40){$E_8$}
\put(220, 40){\circle*{5}}
\put(218, 45){1}
\put(220, 40){\line(1,0){20}}
\put(240, 40){\circle*{5}}
\put(238, 45){3}
\put(240, 40){\line(1,0){20}}
\put(260, 40){\circle*{5}}
\put(258, 45){4}
\put(260, 40){\line(0,-1){20}}
\put(260, 20){\circle*{5}}
\put(265, 18){2}
\put(260, 40){\line(1,0){20}}
\put(280, 40){\circle*{5}}
\put(278, 45){5}
\put(280, 40){\line(1,0){20}}
\put(300, 40){\circle*{5}}
\put(298, 45){6}
\put(300, 40){\line(1,0){20}}
\put(320, 40){\circle*{5}}
\put(318, 45){7}
\put(320, 40){\line(1,0){20}}
\put(340, 40){\circle*{5}}
\put(338, 45){8}

\put( 10, 80){$G_2$}
\put( 40, 80){\circle*{5}}
\put( 38, 85){1}
\put( 40, 78){\line(1,0){20}}
\put( 40, 80){\line(1,0){20}}
\put( 40, 82){\line(1,0){20}}
\put( 46, 78){$>$}
\put( 60, 80){\circle*{5}}
\put( 58, 85){2}

\put(100, 80){$F_4$}
\put(130, 80){\circle*{5}}
\put(128, 85){1}
\put(130, 80){\line(1,0){20}}
\put(150, 80){\circle*{5}}
\put(148, 85){2}
\put(150, 78){\line(1,0){20}}
\put(150, 81){\line(1,0){20}}
\put(156, 77){$>$}
\put(170, 80){\circle*{5}}
\put(168, 85){3}
\put(170, 80){\line(1,0){20}}
\put(190, 80){\circle*{5}}
\put(188, 85){4}

\put(230, 80){$E_6$}
\put(260, 80){\circle*{5}}
\put(258, 85){1}
\put(260, 80){\line(1,0){20}}
\put(280, 80){\circle*{5}}
\put(278, 85){3}
\put(280, 80){\line(1,0){20}}
\put(300, 80){\circle*{5}}
\put(298, 85){4}
\put(300, 80){\line(0,-1){20}}
\put(300, 60){\circle*{5}}
\put(305, 58){2}
\put(300, 80){\line(1,0){20}}
\put(320, 80){\circle*{5}}
\put(318, 85){5}
\put(320, 80){\line(1,0){20}}
\put(340, 80){\circle*{5}}
\put(338, 85){6}

\put( 10,130){$D_4$}
\put( 40,150){\circle*{5}}
\put( 45,150){1}
\put( 40,150){\line(1,-1){21}}
\put( 40,110){\circle*{5}}
\put( 45,108){2}
\put( 40,110){\line(1,1){21}}
\put( 60,130){\circle*{5}}
\put( 58,135){3}
\put( 60,130){\line(1,0){20}}
\put( 80,130){\circle*{5}}
\put( 78,135){4}

\put(10,170){$B_2$}
\put(40,170){\circle*{5}}
\put(38,175){1}
\put(40,168){\line(1,0){20}}
\put(40,171){\line(1,0){20}}
\put(46,167){$<$}
\put(60,170){\circle*{5}}
\put(58,175){2}
\end{picture}
\setlength{\unitlength}{\tmplen}
}
\end{center}

\subsection{$G_2(q)$}

\[
\begin{array}{lll}
\textrm{class}&\textrm{representative}&\textrm{elementary divisors}\\
\hline
A_0 & - & \Phi_{1},\Phi_{1}\\
\tilde{A}_1 & 2 & \Phi_{1}\Phi_{2}\\
A_1 & 1 & \Phi_{1}\Phi_{2}\\
G_2 & 12 & \Phi_{6}\\
A_2 & 1212 & \Phi_{3}\\
A_1+\tilde{A}_1 & 121212 & \Phi_{2},\Phi_{2}\\
\end{array}
\]

\subsection{${}^3D_4(q)$}

\[
\begin{array}{lll}
\textrm{class}&\textrm{representative}&\textrm{elementary divisors}\\
\hline
C_3 & \begin{minipage}[t]{8em}1\-\\\end{minipage} & \Phi_{1}\Phi_{2}\Phi_{6}\\
\tilde{A}_2 & \begin{minipage}[t]{8em}-\-\\\end{minipage} & \Phi_{1},\Phi_{1}\Phi_{3}\\
C_3+A_1 & \begin{minipage}[t]{8em}1\-2\-3\-1\-2\-3\-\\\end{minipage} & \Phi_{2},\Phi_{2}\Phi_{6}\\
\tilde{A}_2+A_1 & \begin{minipage}[t]{8em}3\-\\\end{minipage} & \Phi_{1}\Phi_{2}\Phi_{3}\\
F_4 & \begin{minipage}[t]{8em}1\-3\-\\\end{minipage} & \Phi_{12}\\
\tilde{A}_2+A_2 & \begin{minipage}[t]{8em}1\-2\-3\-1\-2\-4\-3\-2\-\\\end{minipage} & \Phi_{3},\Phi_{3}\\
F_4(a_1) & \begin{minipage}[t]{8em}1\-2\-3\-2\-\\\end{minipage} & \Phi_{6},\Phi_{6}\\
\end{array}
\]

\newpage
\mbox{}\\[-20mm]
\subsection{$F_4(q)$}
\nopagebreak
\[
\begin{array}{lll}
\textrm{class}&\textrm{representative}&\textrm{elementary divisors}\\
\hline
A_0 & \begin{minipage}[t]{13em}-\-\\\end{minipage} & \Phi_{1},\Phi_{1},\Phi_{1},\Phi_{1}\\
4A_1 & \begin{minipage}[t]{13em}1\-2\-1\-3\-2\-1\-3\-2\-3\-4\-3\-2\-1\-3\-2\-3\-4\-3\-2\-1\-3\-2\-3\-4\-\\\end{minipage} & \Phi_{2},\Phi_{2},\Phi_{2},\Phi_{2}\\
2A_1 & \begin{minipage}[t]{13em}2\-3\-2\-3\-\\\end{minipage} & \left\{\begin{array}{ll}
\Phi_{1}\Phi_{2},\Phi_{1}\Phi_{2} & \textrm{ if } q \equiv 0 \mod 2\\
2,1/2 \Phi_{1}\Phi_{2},\Phi_{1}\Phi_{2} & \textrm{ if } q \equiv 1 \mod 2\\
\end{array}\right.
\\
 & & \textrm{uniform: }\Phi_{1},\Phi_{2},\Phi_{1}\Phi_{2}\\
A_2 & \begin{minipage}[t]{13em}2\-1\-\\\end{minipage} & \Phi_{1},\Phi_{1}\Phi_{3}\\
D_4 & \begin{minipage}[t]{13em}2\-3\-2\-3\-4\-3\-2\-1\-3\-4\-\\\end{minipage} & \Phi_{2},\Phi_{2}\Phi_{6}\\
D_4(a_1) & \begin{minipage}[t]{13em}3\-2\-4\-3\-2\-1\-3\-2\-4\-3\-2\-1\-\\\end{minipage} & \Phi_{4},\Phi_{4}\\
\tilde{A}_2 & \begin{minipage}[t]{13em}4\-3\-\\\end{minipage} & \Phi_{1},\Phi_{1}\Phi_{3}\\
C_3+A_1 & \begin{minipage}[t]{13em}1\-2\-1\-4\-3\-2\-1\-3\-2\-3\-\\\end{minipage} & \Phi_{2},\Phi_{2}\Phi_{6}\\
A_2+\tilde{A}_2 & \begin{minipage}[t]{13em}3\-2\-1\-4\-3\-2\-1\-3\-2\-3\-4\-3\-2\-1\-3\-2\-\\\end{minipage} & \Phi_{3},\Phi_{3}\\
F_4(a_1) & \begin{minipage}[t]{13em}3\-2\-4\-3\-2\-1\-3\-2\-\\\end{minipage} & \Phi_{6},\Phi_{6}\\
F_4 & \begin{minipage}[t]{13em}4\-3\-2\-1\-\\\end{minipage} & \Phi_{12}\\
A_1 & \begin{minipage}[t]{13em}1\-\\\end{minipage} & \Phi_{1},\Phi_{1},\Phi_{1}\Phi_{2}\\
3A_1 & \begin{minipage}[t]{13em}2\-3\-2\-3\-4\-3\-2\-3\-4\-\\\end{minipage} & \Phi_{2},\Phi_{2},\Phi_{1}\Phi_{2}\\
\tilde{A}_2+A_1 & \begin{minipage}[t]{13em}1\-4\-3\-\\\end{minipage} & \Phi_{1}\Phi_{2}\Phi_{3}\\
C_3 & \begin{minipage}[t]{13em}4\-3\-2\-\\\end{minipage} & \Phi_{1}\Phi_{2}\Phi_{6}\\
A_3 & \begin{minipage}[t]{13em}2\-3\-2\-1\-3\-\\\end{minipage} & \left\{\begin{array}{ll}
\Phi_{1}\Phi_{2}\Phi_{4} & \textrm{ if } q \equiv 0 \mod 2\\
2,1/2 \Phi_{1}\Phi_{2}\Phi_{4} & \textrm{ if } q \equiv 1 \mod 2\\
\end{array}\right.
\\
 & & \textrm{uniform: }\Phi_{4},\Phi_{1}\Phi_{2}\\
\tilde{A}_1 & \begin{minipage}[t]{13em}3\-\\\end{minipage} & \Phi_{1},\Phi_{1},\Phi_{1}\Phi_{2}\\
2A_1+\tilde{A}_1 & \begin{minipage}[t]{13em}1\-2\-1\-3\-2\-1\-3\-2\-3\-\\\end{minipage} & \Phi_{2},\Phi_{2},\Phi_{1}\Phi_{2}\\
A_2+\tilde{A}_1 & \begin{minipage}[t]{13em}2\-1\-4\-\\\end{minipage} & \Phi_{1}\Phi_{2}\Phi_{3}\\
B_3 & \begin{minipage}[t]{13em}3\-2\-1\-\\\end{minipage} & \Phi_{1}\Phi_{2}\Phi_{6}\\
B_2+A_1 & \begin{minipage}[t]{13em}2\-4\-3\-2\-3\-\\\end{minipage} & \left\{\begin{array}{ll}
\Phi_{1}\Phi_{2}\Phi_{4} & \textrm{ if } q \equiv 0 \mod 2\\
2,1/2 \Phi_{1}\Phi_{2}\Phi_{4} & \textrm{ if } q \equiv 1 \mod 2\\
\end{array}\right.
\\
 & & \textrm{uniform: }\Phi_{4},\Phi_{1}\Phi_{2}\\
A_1+\tilde{A}_1 & \begin{minipage}[t]{13em}1\-3\-\\\end{minipage} & \Phi_{1}\Phi_{2},\Phi_{1}\Phi_{2}\\
B_2 & \begin{minipage}[t]{13em}3\-2\-\\\end{minipage} & \Phi_{1},\Phi_{1}\Phi_{4}\\
A_3+\tilde{A}_1 & \begin{minipage}[t]{13em}2\-3\-2\-3\-4\-3\-2\-1\-3\-2\-4\-3\-2\-1\-\\\end{minipage} & \Phi_{2},\Phi_{2}\Phi_{4}\\
B_4 & \begin{minipage}[t]{13em}2\-4\-3\-2\-1\-3\-\\\end{minipage} & \Phi_{8}\\
\end{array}
\]

\newpage
\mbox{}\\[-20mm]
\subsection{$E_6(q)$}

This holds for the groups of simply-connected type and of adjoint type.

\[
\begin{array}{lll}
\textrm{class}&\textrm{representative}&\textrm{elementary divisors}\\
\hline
A_0 & \begin{minipage}[t]{8em}-\-\\\end{minipage} & \Phi_{1},\Phi_{1},\Phi_{1},\Phi_{1},\Phi_{1},\Phi_{1}\\
4A_1 & \begin{minipage}[t]{8em}3\-4\-3\-2\-4\-3\-5\-4\-3\-2\-4\-5\-\\\end{minipage} & \Phi_{2},\Phi_{2},\Phi_{1}\Phi_{2},\Phi_{1}\Phi_{2}\\
2A_1 & \begin{minipage}[t]{8em}1\-4\-\\\end{minipage} & \Phi_{1},\Phi_{1},\Phi_{1}\Phi_{2},\Phi_{1}\Phi_{2}\\
3A_2 & \begin{minipage}[t]{8em}1\-3\-1\-4\-3\-1\-2\-4\-5\-4\-3\-1\-2\-4\-3\-5\-6\-5\-4\-3\-2\-4\-5\-6\-\\\end{minipage} & \Phi_{3},\Phi_{3},\Phi_{3}\\
A_2 & \begin{minipage}[t]{8em}1\-3\-\\\end{minipage} & \Phi_{1},\Phi_{1},\Phi_{1},\Phi_{1}\Phi_{3}\\
2A_2 & \begin{minipage}[t]{8em}1\-3\-5\-6\-\\\end{minipage} & \left\{\begin{array}{ll}
\Phi_{1}\Phi_{3},\Phi_{1}\Phi_{3} & \textrm{ if } q \equiv 0,2 \mod 3\\
3,1/3 \Phi_{1}\Phi_{3},\Phi_{1}\Phi_{3} & \textrm{ if } q \equiv 1 \mod 3\\
\end{array}\right.
\\
 & & \textrm{uniform: }\Phi_{1},\Phi_{3},\Phi_{1}\Phi_{3}\\
D_4(a_1) & \begin{minipage}[t]{8em}3\-4\-3\-2\-4\-5\-\\\end{minipage} & \Phi_{1}\Phi_{4},\Phi_{1}\Phi_{4}\\
A_3+A_1 & \begin{minipage}[t]{8em}1\-4\-3\-6\-\\\end{minipage} & \Phi_{1}\Phi_{2},\Phi_{1}\Phi_{2}\Phi_{4}\\
A_4 & \begin{minipage}[t]{8em}1\-4\-3\-2\-\\\end{minipage} & \Phi_{1},\Phi_{1}\Phi_{5}\\
E_6(a_2) & \begin{minipage}[t]{8em}1\-2\-3\-1\-5\-4\-6\-5\-4\-2\-3\-4\-\\\end{minipage} & \Phi_{6},\Phi_{3}\Phi_{6}\\
D_4 & \begin{minipage}[t]{8em}3\-4\-2\-5\-\\\end{minipage} & \Phi_{1}\Phi_{2},\Phi_{1}\Phi_{2}\Phi_{6}\\
A_5+A_1 & \begin{minipage}[t]{8em}1\-2\-3\-4\-2\-3\-4\-6\-5\-4\-2\-3\-4\-5\-\\\end{minipage} & \Phi_{2},\Phi_{2}\Phi_{3}\Phi_{6}\\
A_2+2A_1 & \begin{minipage}[t]{8em}1\-3\-2\-5\-\\\end{minipage} & \Phi_{1}\Phi_{2},\Phi_{1}\Phi_{2}\Phi_{3}\\
E_6(a_1) & \begin{minipage}[t]{8em}1\-3\-4\-3\-2\-4\-5\-6\-\\\end{minipage} & \Phi_{9}\\
E_6 & \begin{minipage}[t]{8em}1\-4\-6\-2\-3\-5\-\\\end{minipage} & \Phi_{3}\Phi_{12}\\
A_1 & \begin{minipage}[t]{8em}1\-\\\end{minipage} & \Phi_{1},\Phi_{1},\Phi_{1},\Phi_{1},\Phi_{1}\Phi_{2}\\
3A_1 & \begin{minipage}[t]{8em}1\-4\-6\-\\\end{minipage} & \Phi_{1}\Phi_{2},\Phi_{1}\Phi_{2},\Phi_{1}\Phi_{2}\\
A_3+2A_1 & \begin{minipage}[t]{8em}1\-3\-4\-3\-2\-4\-3\-5\-4\-3\-2\-4\-5\-\\\end{minipage} & \Phi_{2},\Phi_{2},\Phi_{1}\Phi_{2}\Phi_{4}\\
A_3 & \begin{minipage}[t]{8em}1\-4\-3\-\\\end{minipage} & \Phi_{1},\Phi_{1},\Phi_{1}\Phi_{2}\Phi_{4}\\
A_2+A_1 & \begin{minipage}[t]{8em}1\-3\-2\-\\\end{minipage} & \Phi_{1},\Phi_{1},\Phi_{1}\Phi_{2}\Phi_{3}\\
2A_2+A_1 & \begin{minipage}[t]{8em}1\-3\-2\-5\-6\-\\\end{minipage} & \Phi_{3},\Phi_{1}\Phi_{2}\Phi_{3}\\
A_5 & \begin{minipage}[t]{8em}1\-4\-6\-3\-5\-\\\end{minipage} & \left\{\begin{array}{ll}
\Phi_{1}\Phi_{2}\Phi_{3}\Phi_{6} & \textrm{ if } q \equiv 0,2 \mod 3\\
3,1/3 \Phi_{1}\Phi_{2}\Phi_{3}\Phi_{6} & \textrm{ if } q \equiv 1 \mod 3\\
\end{array}\right.
\\
 & & \textrm{uniform: }\Phi_{3},\Phi_{1}\Phi_{2}\Phi_{6}\\
D_5 & \begin{minipage}[t]{8em}1\-3\-4\-2\-5\-\\\end{minipage} & \Phi_{1}\Phi_{2}\Phi_{8}\\
A_4+A_1 & \begin{minipage}[t]{8em}1\-4\-3\-2\-6\-\\\end{minipage} & \Phi_{1}\Phi_{2}\Phi_{5}\\
D_5(a_1) & \begin{minipage}[t]{8em}1\-4\-2\-5\-4\-2\-3\-\\\end{minipage} & \Phi_{1}\Phi_{2}\Phi_{4}\Phi_{6}\\
\end{array}
\]

And  here  is the  structure  of  the maximal  tori  in  the simple  group
$(E_6)_{sc}(q)/Z$ in the case $q \equiv 1 \mod 3$ (then $|Z| = 3$).

\[
\begin{array}{lll}
\textrm{class}&\textrm{representative}&\textrm{elementary divisors}\\
\hline
A_0 & \begin{minipage}[t]{8em}-\-\\\end{minipage} & 1/3 \Phi_{1},\Phi_{1},\Phi_{1},\Phi_{1},\Phi_{1},\Phi_{1}\\
4A_1 & \begin{minipage}[t]{8em}3\-4\-3\-2\-4\-3\-5\-4\-3\-2\-4\-5\-\\\end{minipage} & \Phi_{2},\Phi_{2},1/3 \Phi_{1}\Phi_{2},\Phi_{1}\Phi_{2}\\
2A_1 & \begin{minipage}[t]{8em}1\-4\-\\\end{minipage} & 1/3 \Phi_{1},\Phi_{1},\Phi_{1}\Phi_{2},\Phi_{1}\Phi_{2}\\
3A_2 & \begin{minipage}[t]{8em}1\-3\-1\-4\-3\-1\-2\-4\-5\-4\-3\-1\-2\-4\-3\-5\-6\-5\-4\-3\-2\-4\-5\-6\-\\\end{minipage} & 1/3 \Phi_{3},\Phi_{3},\Phi_{3}\\
A_2 & \begin{minipage}[t]{8em}1\-3\-\\\end{minipage} & 1/3 \Phi_{1},\Phi_{1},\Phi_{1},\Phi_{1}\Phi_{3}\\
2A_2 & \begin{minipage}[t]{8em}1\-3\-5\-6\-\\\end{minipage} & 1/3 \Phi_{1}\Phi_{3},\Phi_{1}\Phi_{3}\\
D_4(a_1) & \begin{minipage}[t]{8em}3\-4\-3\-2\-4\-5\-\\\end{minipage} & 1/3 \Phi_{1}\Phi_{4},\Phi_{1}\Phi_{4}\\
A_3+A_1 & \begin{minipage}[t]{8em}1\-4\-3\-6\-\\\end{minipage} & 1/3 \Phi_{1}\Phi_{2},\Phi_{1}\Phi_{2}\Phi_{4}\\
A_4 & \begin{minipage}[t]{8em}1\-4\-3\-2\-\\\end{minipage} & 1/3 \Phi_{1},\Phi_{1}\Phi_{5}\\
E_6(a_2) & \begin{minipage}[t]{8em}1\-2\-3\-1\-5\-4\-6\-5\-4\-2\-3\-4\-\\\end{minipage} & \Phi_{6},1/3 \Phi_{3}\Phi_{6}\\
D_4 & \begin{minipage}[t]{8em}3\-4\-2\-5\-\\\end{minipage} & 1/3 \Phi_{1}\Phi_{2},\Phi_{1}\Phi_{2}\Phi_{6}\\
A_5+A_1 & \begin{minipage}[t]{8em}1\-2\-3\-4\-2\-3\-4\-6\-5\-4\-2\-3\-4\-5\-\\\end{minipage} & \Phi_{2},1/3 \Phi_{2}\Phi_{3}\Phi_{6}\\
A_2+2A_1 & \begin{minipage}[t]{8em}1\-3\-2\-5\-\\\end{minipage} & 1/3 \Phi_{1}\Phi_{2},\Phi_{1}\Phi_{2}\Phi_{3}\\
E_6(a_1) & \begin{minipage}[t]{8em}1\-3\-4\-3\-2\-4\-5\-6\-\\\end{minipage} & 1/3 \Phi_{9}\\
E_6 & \begin{minipage}[t]{8em}1\-4\-6\-2\-3\-5\-\\\end{minipage} & 1/3 \Phi_{3}\Phi_{12}\\
A_1 & \begin{minipage}[t]{8em}1\-\\\end{minipage} & 1/3 \Phi_{1},\Phi_{1},\Phi_{1},\Phi_{1},\Phi_{1}\Phi_{2}\\
3A_1 & \begin{minipage}[t]{8em}1\-4\-6\-\\\end{minipage} & 1/3 \Phi_{1}\Phi_{2},\Phi_{1}\Phi_{2},\Phi_{1}\Phi_{2}\\
A_3+2A_1 & \begin{minipage}[t]{8em}1\-3\-4\-3\-2\-4\-3\-5\-4\-3\-2\-4\-5\-\\\end{minipage} & \Phi_{2},\Phi_{2},1/3 \Phi_{1}\Phi_{2}\Phi_{4}\\
A_3 & \begin{minipage}[t]{8em}1\-4\-3\-\\\end{minipage} & 1/3 \Phi_{1},\Phi_{1},\Phi_{1}\Phi_{2}\Phi_{4}\\
A_2+A_1 & \begin{minipage}[t]{8em}1\-3\-2\-\\\end{minipage} & 1/3 \Phi_{1},\Phi_{1},\Phi_{1}\Phi_{2}\Phi_{3}\\
2A_2+A_1 & \begin{minipage}[t]{8em}1\-3\-2\-5\-6\-\\\end{minipage} & 1/3 \Phi_{3},\Phi_{1}\Phi_{2}\Phi_{3}\\
A_5 & \begin{minipage}[t]{8em}1\-4\-6\-3\-5\-\\\end{minipage} & 1/3 \Phi_{1}\Phi_{2}\Phi_{3}\Phi_{6}\\
D_5 & \begin{minipage}[t]{8em}1\-3\-4\-2\-5\-\\\end{minipage} & 1/3 \Phi_{1}\Phi_{2}\Phi_{8}\\
A_4+A_1 & \begin{minipage}[t]{8em}1\-4\-3\-2\-6\-\\\end{minipage} & 1/3 \Phi_{1}\Phi_{2}\Phi_{5}\\
D_5(a_1) & \begin{minipage}[t]{8em}1\-4\-2\-5\-4\-2\-3\-\\\end{minipage} & 1/3 \Phi_{1}\Phi_{2}\Phi_{4}\Phi_{6}\\
\end{array}
\]

\newpage
\mbox{}\\[-20mm]
\subsection{${}^2E_6(q)$}

This holds for the groups of simply-connected type and of adjoint type.

\[
\begin{array}{lll}
\textrm{class}&\textrm{representative}&\textrm{elementary divisors}\\
\hline
A_0 & \begin{minipage}[t]{10em}1\-2\-3\-1\-4\-2\-3\-1\-4\-3\-5\-4\-2\-3\-1\-4\-3\-5\-4\-2\-6\-5\-4\-2\-3\-1\-4\-3\-5\-4\-2\-6\-5\-4\-3\-1\-\\\end{minipage} & \Phi_{2},\Phi_{2},\Phi_{2},\Phi_{2},\Phi_{2},\Phi_{2}\\
4A_1 & \begin{minipage}[t]{10em}-\-\\\end{minipage} & \Phi_{1},\Phi_{1},\Phi_{1}\Phi_{2},\Phi_{1}\Phi_{2}\\
2A_1 & \begin{minipage}[t]{10em}3\-4\-3\-5\-4\-3\-\\\end{minipage} & \Phi_{2},\Phi_{2},\Phi_{1}\Phi_{2},\Phi_{1}\Phi_{2}\\
3A_2 & \begin{minipage}[t]{10em}1\-2\-4\-3\-1\-5\-4\-3\-6\-5\-4\-3\-\\\end{minipage} & \Phi_{6},\Phi_{6},\Phi_{6}\\
A_2 & \begin{minipage}[t]{10em}1\-2\-3\-1\-4\-3\-1\-5\-4\-3\-1\-6\-5\-4\-3\-1\-\\\end{minipage} & \Phi_{2},\Phi_{2},\Phi_{2},\Phi_{2}\Phi_{6}\\
2A_2 & \begin{minipage}[t]{10em}2\-3\-4\-2\-3\-4\-6\-5\-4\-2\-3\-4\-5\-6\-\\\end{minipage} & \left\{\begin{array}{ll}
\Phi_{2}\Phi_{6},\Phi_{2}\Phi_{6} & \textrm{ if } q \equiv 0,1 \mod 3\\
3,1/3 \Phi_{2}\Phi_{6},\Phi_{2}\Phi_{6} & \textrm{ if } q \equiv 2 \mod 3\\
\end{array}\right.
\\
 & & \textrm{uniform: }\Phi_{2},\Phi_{6},\Phi_{2}\Phi_{6}\\
D_4(a_1) & \begin{minipage}[t]{10em}1\-4\-2\-3\-1\-4\-3\-5\-4\-2\-3\-1\-4\-6\-5\-4\-3\-1\-\\\end{minipage} & \Phi_{2}\Phi_{4},\Phi_{2}\Phi_{4}\\
A_3+A_1 & \begin{minipage}[t]{10em}1\-2\-\\\end{minipage} & \Phi_{1}\Phi_{2},\Phi_{1}\Phi_{2}\Phi_{4}\\
A_4 & \begin{minipage}[t]{10em}4\-5\-4\-2\-3\-1\-4\-5\-\\\end{minipage} & \Phi_{2},\Phi_{2}\Phi_{10}\\
E_6(a_2) & \begin{minipage}[t]{10em}4\-2\-5\-4\-2\-3\-4\-5\-6\-5\-4\-2\-3\-4\-5\-6\-\\\end{minipage} & \Phi_{3},\Phi_{3}\Phi_{6}\\
D_4 & \begin{minipage}[t]{10em}2\-4\-\\\end{minipage} & \Phi_{1}\Phi_{2},\Phi_{1}\Phi_{2}\Phi_{3}\\
A_5+A_1 & \begin{minipage}[t]{10em}1\-5\-\\\end{minipage} & \Phi_{1},\Phi_{1}\Phi_{3}\Phi_{6}\\
A_2+2A_1 & \begin{minipage}[t]{10em}5\-4\-\\\end{minipage} & \Phi_{1}\Phi_{2},\Phi_{1}\Phi_{2}\Phi_{6}\\
E_6(a_1) & \begin{minipage}[t]{10em}1\-2\-5\-4\-\\\end{minipage} & \Phi_{18}\\
E_6 & \begin{minipage}[t]{10em}1\-2\-3\-1\-4\-3\-\\\end{minipage} & \Phi_{6}\Phi_{12}\\
A_1 & \begin{minipage}[t]{10em}1\-3\-1\-4\-3\-1\-5\-4\-3\-1\-6\-5\-4\-3\-1\-\\\end{minipage} & \Phi_{2},\Phi_{2},\Phi_{2},\Phi_{2},\Phi_{1}\Phi_{2}\\
3A_1 & \begin{minipage}[t]{10em}2\-\\\end{minipage} & \Phi_{1}\Phi_{2},\Phi_{1}\Phi_{2},\Phi_{1}\Phi_{2}\\
A_3+2A_1 & \begin{minipage}[t]{10em}1\-\\\end{minipage} & \Phi_{1},\Phi_{1},\Phi_{1}\Phi_{2}\Phi_{4}\\
A_3 & \begin{minipage}[t]{10em}2\-3\-4\-3\-5\-4\-3\-\\\end{minipage} & \Phi_{2},\Phi_{2},\Phi_{1}\Phi_{2}\Phi_{4}\\
A_2+A_1 & \begin{minipage}[t]{10em}1\-3\-4\-3\-5\-4\-3\-\\\end{minipage} & \Phi_{2},\Phi_{2},\Phi_{1}\Phi_{2}\Phi_{6}\\
2A_2+A_1 & \begin{minipage}[t]{10em}1\-3\-1\-4\-3\-\\\end{minipage} & \Phi_{6},\Phi_{1}\Phi_{2}\Phi_{6}\\
A_5 & \begin{minipage}[t]{10em}1\-2\-5\-\\\end{minipage} & \left\{\begin{array}{ll}
\Phi_{1}\Phi_{2}\Phi_{3}\Phi_{6} & \textrm{ if } q \equiv 0,1 \mod 3\\
3,1/3 \Phi_{1}\Phi_{2}\Phi_{3}\Phi_{6} & \textrm{ if } q \equiv 2 \mod 3\\
\end{array}\right.
\\
 & & \textrm{uniform: }\Phi_{6},\Phi_{1}\Phi_{2}\Phi_{3}\\
D_5 & \begin{minipage}[t]{10em}2\-5\-4\-\\\end{minipage} & \Phi_{1}\Phi_{2}\Phi_{8}\\
A_4+A_1 & \begin{minipage}[t]{10em}1\-5\-4\-\\\end{minipage} & \Phi_{1}\Phi_{2}\Phi_{10}\\
D_5(a_1) & \begin{minipage}[t]{10em}1\-2\-4\-\\\end{minipage} & \Phi_{1}\Phi_{2}\Phi_{3}\Phi_{4}\\
\end{array}
\]

And  here  is the  structure  of  the maximal  tori  in  the simple  group
$({}^2E_6)_{sc}(q)/Z$ in the case $q \equiv 2 \mod 3$ (then $|Z| = 3$).

\[
\begin{array}{lll}
\textrm{class}&\textrm{representative}&\textrm{elementary divisors}\\
\hline
A_0 & \begin{minipage}[t]{10em}1\-2\-3\-1\-4\-2\-3\-1\-4\-3\-5\-4\-2\-3\-1\-4\-3\-5\-4\-2\-6\-5\-4\-2\-3\-1\-4\-3\-5\-4\-2\-6\-5\-4\-3\-1\-\\\end{minipage} & 1/3 \Phi_{2},\Phi_{2},\Phi_{2},\Phi_{2},\Phi_{2},\Phi_{2}\\
4A_1 & \begin{minipage}[t]{10em}-\-\\\end{minipage} & \Phi_{1},\Phi_{1},1/3 \Phi_{1}\Phi_{2},\Phi_{1}\Phi_{2}\\
2A_1 & \begin{minipage}[t]{10em}3\-4\-3\-5\-4\-3\-\\\end{minipage} & 1/3 \Phi_{2},\Phi_{2},\Phi_{1}\Phi_{2},\Phi_{1}\Phi_{2}\\
3A_2 & \begin{minipage}[t]{10em}1\-2\-4\-3\-1\-5\-4\-3\-6\-5\-4\-3\-\\\end{minipage} & 1/3 \Phi_{6},\Phi_{6},\Phi_{6}\\
A_2 & \begin{minipage}[t]{10em}1\-2\-3\-1\-4\-3\-1\-5\-4\-3\-1\-6\-5\-4\-3\-1\-\\\end{minipage} & 1/3 \Phi_{2},\Phi_{2},\Phi_{2},\Phi_{2}\Phi_{6}\\
2A_2 & \begin{minipage}[t]{10em}2\-3\-4\-2\-3\-4\-6\-5\-4\-2\-3\-4\-5\-6\-\\\end{minipage} & 1/3 \Phi_{2}\Phi_{6},\Phi_{2}\Phi_{6}\\
D_4(a_1) & \begin{minipage}[t]{10em}1\-4\-2\-3\-1\-4\-3\-5\-4\-2\-3\-1\-4\-6\-5\-4\-3\-1\-\\\end{minipage} & 1/3 \Phi_{2}\Phi_{4},\Phi_{2}\Phi_{4}\\
A_3+A_1 & \begin{minipage}[t]{10em}1\-2\-\\\end{minipage} & 1/3 \Phi_{1}\Phi_{2},\Phi_{1}\Phi_{2}\Phi_{4}\\
A_4 & \begin{minipage}[t]{10em}4\-5\-4\-2\-3\-1\-4\-5\-\\\end{minipage} & 1/3 \Phi_{2},\Phi_{2}\Phi_{10}\\
E_6(a_2) & \begin{minipage}[t]{10em}4\-2\-5\-4\-2\-3\-4\-5\-6\-5\-4\-2\-3\-4\-5\-6\-\\\end{minipage} & \Phi_{3},1/3 \Phi_{3}\Phi_{6}\\
D_4 & \begin{minipage}[t]{10em}2\-4\-\\\end{minipage} & 1/3 \Phi_{1}\Phi_{2},\Phi_{1}\Phi_{2}\Phi_{3}\\
A_5+A_1 & \begin{minipage}[t]{10em}1\-5\-\\\end{minipage} & \Phi_{1},1/3 \Phi_{1}\Phi_{3}\Phi_{6}\\
A_2+2A_1 & \begin{minipage}[t]{10em}5\-4\-\\\end{minipage} & 1/3 \Phi_{1}\Phi_{2},\Phi_{1}\Phi_{2}\Phi_{6}\\
E_6(a_1) & \begin{minipage}[t]{10em}1\-2\-5\-4\-\\\end{minipage} & 1/3 \Phi_{18}\\
E_6 & \begin{minipage}[t]{10em}1\-2\-3\-1\-4\-3\-\\\end{minipage} & 1/3 \Phi_{6}\Phi_{12}\\
A_1 & \begin{minipage}[t]{10em}1\-3\-1\-4\-3\-1\-5\-4\-3\-1\-6\-5\-4\-3\-1\-\\\end{minipage} & 1/3 \Phi_{2},\Phi_{2},\Phi_{2},\Phi_{2},\Phi_{1}\Phi_{2}\\
3A_1 & \begin{minipage}[t]{10em}2\-\\\end{minipage} & 1/3 \Phi_{1}\Phi_{2},\Phi_{1}\Phi_{2},\Phi_{1}\Phi_{2}\\
A_3+2A_1 & \begin{minipage}[t]{10em}1\-\\\end{minipage} & \Phi_{1},\Phi_{1},1/3 \Phi_{1}\Phi_{2}\Phi_{4}\\
A_3 & \begin{minipage}[t]{10em}2\-3\-4\-3\-5\-4\-3\-\\\end{minipage} & 1/3 \Phi_{2},\Phi_{2},\Phi_{1}\Phi_{2}\Phi_{4}\\
A_2+A_1 & \begin{minipage}[t]{10em}1\-3\-4\-3\-5\-4\-3\-\\\end{minipage} & 1/3 \Phi_{2},\Phi_{2},\Phi_{1}\Phi_{2}\Phi_{6}\\
2A_2+A_1 & \begin{minipage}[t]{10em}1\-3\-1\-4\-3\-\\\end{minipage} & 1/3 \Phi_{6},\Phi_{1}\Phi_{2}\Phi_{6}\\
A_5 & \begin{minipage}[t]{10em}1\-2\-5\-\\\end{minipage} & 1/3 \Phi_{1}\Phi_{2}\Phi_{3}\Phi_{6}\\
D_5 & \begin{minipage}[t]{10em}2\-5\-4\-\\\end{minipage} & 1/3 \Phi_{1}\Phi_{2}\Phi_{8}\\
A_4+A_1 & \begin{minipage}[t]{10em}1\-5\-4\-\\\end{minipage} & 1/3 \Phi_{1}\Phi_{2}\Phi_{10}\\
D_5(a_1) & \begin{minipage}[t]{10em}1\-2\-4\-\\\end{minipage} & 1/3 \Phi_{1}\Phi_{2}\Phi_{3}\Phi_{4}\\
\end{array}
\]

\newpage
\subsection{$E_7(q)$}

This holds for the groups of simply-connected type and of adjoint type.

\[
\begin{array}{lll}
\textrm{class}&\textrm{representative}&\textrm{elementary divisors}\\
\hline
A_0 & \begin{minipage}[t]{10em}-\-\\\end{minipage} & \Phi_{1},\Phi_{1},\Phi_{1},\Phi_{1},\Phi_{1},\Phi_{1},\Phi_{1}\\
6A_1 & \begin{minipage}[t]{10em}7\-6\-7\-5\-6\-7\-4\-5\-6\-7\-2\-4\-5\-6\-7\-3\-4\-5\-6\-7\-2\-4\-5\-6\-3\-4\-5\-2\-4\-3\-\\\end{minipage} & \Phi_{2},\Phi_{2},\Phi_{2},\Phi_{2},\Phi_{2},\Phi_{1}\Phi_{2}\\
4A_1'' & \begin{minipage}[t]{10em}5\-4\-5\-2\-4\-5\-3\-4\-5\-2\-4\-3\-\\\end{minipage} & \left\{\begin{array}{ll}
\Phi_{2},\Phi_{1}\Phi_{2},\Phi_{1}\Phi_{2},\Phi_{1}\Phi_{2} & \textrm{ if } q \equiv 0 \mod 2\\
2,\Phi_{2},1/2 \Phi_{1}\Phi_{2},\Phi_{1}\Phi_{2},\Phi_{1}\Phi_{2} & \textrm{ if } q \equiv 1 \mod 2\\
\end{array}\right.
\\
 & & \textrm{uniform: }\Phi_{1},\Phi_{2},\Phi_{2},\Phi_{1}\Phi_{2},\Phi_{1}\Phi_{2}\\
2A_1 & \begin{minipage}[t]{10em}7\-5\-\\\end{minipage} & \Phi_{1},\Phi_{1},\Phi_{1},\Phi_{1}\Phi_{2},\Phi_{1}\Phi_{2}\\
4A_1' & \begin{minipage}[t]{10em}7\-5\-2\-3\-\\\end{minipage} & \Phi_{2},\Phi_{1}\Phi_{2},\Phi_{1}\Phi_{2},\Phi_{1}\Phi_{2}\\
A_2 & \begin{minipage}[t]{10em}7\-6\-\\\end{minipage} & \Phi_{1},\Phi_{1},\Phi_{1},\Phi_{1},\Phi_{1}\Phi_{3}\\
3A_2 & \begin{minipage}[t]{10em}6\-5\-6\-4\-5\-6\-2\-4\-3\-4\-5\-6\-2\-4\-5\-3\-1\-3\-4\-5\-2\-4\-3\-1\-\\\end{minipage} & \Phi_{3},\Phi_{3},\Phi_{1}\Phi_{3}\\
2A_2 & \begin{minipage}[t]{10em}7\-6\-4\-2\-\\\end{minipage} & \Phi_{1},\Phi_{1}\Phi_{3},\Phi_{1}\Phi_{3}\\
D_4(a_1) & \begin{minipage}[t]{10em}5\-4\-5\-2\-4\-3\-\\\end{minipage} & \Phi_{1},\Phi_{1}\Phi_{4},\Phi_{1}\Phi_{4}\\
A_3+A_1' & \begin{minipage}[t]{10em}7\-5\-6\-2\-\\\end{minipage} & \left\{\begin{array}{ll}
\Phi_{1},\Phi_{1}\Phi_{2},\Phi_{1}\Phi_{2}\Phi_{4} & \textrm{ if } q \equiv 0 \mod 2\\
2,\Phi_{1},1/2 \Phi_{1}\Phi_{2},\Phi_{1}\Phi_{2}\Phi_{4} & \textrm{ if } q \equiv 1 \mod 2\\
\end{array}\right.
\\
 & & \textrm{uniform: }\Phi_{1},\Phi_{2},\Phi_{1},\Phi_{1}\Phi_{2}\Phi_{4}\\
A_3+3A_1 & \begin{minipage}[t]{10em}7\-5\-4\-5\-2\-4\-5\-3\-4\-5\-2\-4\-3\-1\-\\\end{minipage} & \Phi_{2},\Phi_{2},\Phi_{2},\Phi_{1}\Phi_{2}\Phi_{4}\\
D_4(a_1)+2A_1 & \begin{minipage}[t]{10em}7\-6\-7\-5\-6\-4\-5\-2\-4\-5\-3\-4\-5\-2\-4\-3\-\\\end{minipage} & \left\{\begin{array}{ll}
\Phi_{2}\Phi_{4},\Phi_{1}\Phi_{2}\Phi_{4} & \textrm{ if } q \equiv 0 \mod 2\\
2,\Phi_{2}\Phi_{4},1/2 \Phi_{1}\Phi_{2}\Phi_{4} & \textrm{ if } q \equiv 1 \mod 2\\
\end{array}\right.
\\
 & & \textrm{uniform: }\Phi_{4},\Phi_{1}\Phi_{2},\Phi_{2}\Phi_{4}\\
A_3+A_1'' & \begin{minipage}[t]{10em}7\-5\-6\-3\-\\\end{minipage} & \Phi_{1},\Phi_{1}\Phi_{2},\Phi_{1}\Phi_{2}\Phi_{4}\\
A_4 & \begin{minipage}[t]{10em}7\-6\-5\-4\-\\\end{minipage} & \Phi_{1},\Phi_{1},\Phi_{1}\Phi_{5}\\
D_4+2A_1 & \begin{minipage}[t]{10em}7\-6\-5\-4\-5\-2\-4\-5\-3\-4\-5\-2\-4\-3\-\\\end{minipage} & \Phi_{2},\Phi_{2},\Phi_{2},\Phi_{1}\Phi_{2}\Phi_{6}\\
D_4 & \begin{minipage}[t]{10em}5\-4\-2\-3\-\\\end{minipage} & \Phi_{1},\Phi_{1}\Phi_{2},\Phi_{1}\Phi_{2}\Phi_{6}\\
E_6(a_2) & \begin{minipage}[t]{10em}1\-2\-3\-1\-5\-4\-6\-5\-4\-2\-3\-4\-\\\end{minipage} & \Phi_{6},\Phi_{1}\Phi_{3}\Phi_{6}\\
A_2+2A_1 & \begin{minipage}[t]{10em}7\-6\-4\-1\-\\\end{minipage} & \Phi_{1},\Phi_{1}\Phi_{2},\Phi_{1}\Phi_{2}\Phi_{3}\\
D_6(a_2) & \begin{minipage}[t]{10em}7\-6\-7\-5\-6\-4\-5\-2\-4\-3\-\\\end{minipage} & \Phi_{2}\Phi_{6},\Phi_{1}\Phi_{2}\Phi_{6}\\
A_5+A_1'' & \begin{minipage}[t]{10em}1\-2\-3\-4\-2\-3\-4\-6\-5\-4\-2\-3\-4\-5\-\\\end{minipage} & \left\{\begin{array}{ll}
\Phi_{2},\Phi_{1}\Phi_{2}\Phi_{3}\Phi_{6} & \textrm{ if } q \equiv 0 \mod 2\\
2,\Phi_{2},1/2 \Phi_{1}\Phi_{2}\Phi_{3}\Phi_{6} & \textrm{ if } q \equiv 1 \mod 2\\
\end{array}\right.
\\
 & & \textrm{uniform: }\Phi_{2},\Phi_{1}\Phi_{3},\Phi_{2}\Phi_{6}\\
A_5+A_1' & \begin{minipage}[t]{10em}7\-5\-2\-6\-4\-1\-\\\end{minipage} & \Phi_{2},\Phi_{1}\Phi_{2}\Phi_{3}\Phi_{6}\\
A_6 & \begin{minipage}[t]{10em}7\-6\-5\-4\-3\-1\-\\\end{minipage} & \Phi_{1}\Phi_{7}\\
\end{array}\]
\[
\begin{array}{lll}
\textrm{}&\textrm{}&\textrm{cont. }E_7(q)\\
\hline
D_5+A_1 & \begin{minipage}[t]{10em}7\-5\-2\-3\-4\-1\-\\\end{minipage} & \Phi_{2},\Phi_{1}\Phi_{2}\Phi_{8}\\
D_6(a_1) & \begin{minipage}[t]{10em}3\-4\-2\-3\-4\-7\-6\-5\-\\\end{minipage} & \left\{\begin{array}{ll}
\Phi_{1}\Phi_{4}\Phi_{8} & \textrm{ if } q \equiv 0 \mod 2\\
2,1/2 \Phi_{1}\Phi_{4}\Phi_{8} & \textrm{ if } q \equiv 1 \mod 2\\
\end{array}\right.
\\
 & & \textrm{uniform: }\Phi_{4},\Phi_{1}\Phi_{8}\\
E_6(a_1) & \begin{minipage}[t]{10em}6\-5\-4\-5\-2\-4\-3\-1\-\\\end{minipage} & \Phi_{1}\Phi_{9}\\
D_6 & \begin{minipage}[t]{10em}7\-6\-5\-4\-2\-3\-\\\end{minipage} & \Phi_{2},\Phi_{1}\Phi_{2}\Phi_{10}\\
A_3+A_2+A_1 & \begin{minipage}[t]{10em}7\-5\-6\-2\-3\-1\-\\\end{minipage} & \Phi_{2},\Phi_{1}\Phi_{2}\Phi_{3}\Phi_{4}\\
D_5(a_1)+A_1 & \begin{minipage}[t]{10em}7\-5\-4\-5\-2\-4\-3\-1\-\\\end{minipage} & \Phi_{2},\Phi_{1}\Phi_{2}\Phi_{4}\Phi_{6}\\
E_6 & \begin{minipage}[t]{10em}6\-4\-1\-5\-3\-2\-\\\end{minipage} & \Phi_{1}\Phi_{3}\Phi_{12}\\
A_4+A_2 & \begin{minipage}[t]{10em}7\-6\-4\-2\-3\-1\-\\\end{minipage} & \Phi_{1}\Phi_{3}\Phi_{5}\\
7A_1 & \begin{minipage}[t]{10em}7\-6\-7\-5\-6\-7\-4\-5\-6\-7\-2\-4\-5\-6\-7\-3\-4\-5\-6\-7\-2\-4\-5\-6\-3\-4\-5\-2\-4\-3\-1\-3\-4\-5\-6\-7\-2\-4\-5\-6\-3\-4\-5\-2\-4\-3\-1\-3\-4\-5\-6\-7\-2\-4\-5\-6\-3\-4\-5\-2\-4\-3\-1\-\\\end{minipage} & \Phi_{2},\Phi_{2},\Phi_{2},\Phi_{2},\Phi_{2},\Phi_{2},\Phi_{2}\\
A_1 & \begin{minipage}[t]{10em}7\-\\\end{minipage} & \Phi_{1},\Phi_{1},\Phi_{1},\Phi_{1},\Phi_{1},\Phi_{1}\Phi_{2}\\
3A_1' & \begin{minipage}[t]{10em}7\-5\-2\-\\\end{minipage} & \left\{\begin{array}{ll}
\Phi_{1},\Phi_{1}\Phi_{2},\Phi_{1}\Phi_{2},\Phi_{1}\Phi_{2} & \textrm{ if } q \equiv 0 \mod 2\\
2,\Phi_{1},1/2 \Phi_{1}\Phi_{2},\Phi_{1}\Phi_{2},\Phi_{1}\Phi_{2} & \textrm{ if } q \equiv 1 \mod 2\\
\end{array}\right.
\\
 & & \textrm{uniform: }\Phi_{1},\Phi_{1},\Phi_{2},\Phi_{1}\Phi_{2},\Phi_{1}\Phi_{2}\\
5A_1 & \begin{minipage}[t]{10em}7\-5\-4\-5\-2\-4\-5\-3\-4\-5\-2\-4\-3\-\\\end{minipage} & \Phi_{2},\Phi_{2},\Phi_{2},\Phi_{1}\Phi_{2},\Phi_{1}\Phi_{2}\\
3A_1'' & \begin{minipage}[t]{10em}7\-5\-3\-\\\end{minipage} & \Phi_{1},\Phi_{1}\Phi_{2},\Phi_{1}\Phi_{2},\Phi_{1}\Phi_{2}\\
D_4+3A_1 & \begin{minipage}[t]{10em}7\-6\-7\-5\-6\-7\-4\-5\-6\-7\-2\-4\-5\-6\-7\-3\-4\-5\-6\-7\-2\-4\-5\-6\-3\-4\-5\-2\-4\-3\-1\-\\\end{minipage} & \Phi_{2},\Phi_{2},\Phi_{2},\Phi_{2},\Phi_{2}\Phi_{6}\\
E_7(a_4) & \begin{minipage}[t]{10em}7\-6\-7\-5\-6\-7\-4\-5\-6\-2\-4\-3\-4\-5\-2\-4\-3\-1\-3\-4\-2\-\\\end{minipage} & \Phi_{6},\Phi_{6},\Phi_{2}\Phi_{6}\\
D_6(a_2)+A_1 & \begin{minipage}[t]{10em}7\-6\-7\-5\-6\-4\-5\-6\-2\-4\-5\-6\-3\-4\-5\-6\-1\-3\-4\-5\-2\-4\-3\-\\\end{minipage} & \Phi_{2},\Phi_{2}\Phi_{6},\Phi_{2}\Phi_{6}\\
2A_3+A_1 & \begin{minipage}[t]{10em}7\-6\-7\-5\-6\-7\-4\-5\-6\-7\-2\-4\-5\-6\-7\-3\-4\-5\-6\-7\-2\-4\-5\-6\-3\-1\-3\-4\-5\-2\-4\-3\-1\-\\\end{minipage} & \Phi_{2},\Phi_{2}\Phi_{4},\Phi_{2}\Phi_{4}\\
A_3+2A_1'' & \begin{minipage}[t]{10em}6\-5\-4\-5\-2\-4\-5\-3\-4\-5\-2\-4\-3\-\\\end{minipage} & \left\{\begin{array}{ll}
\Phi_{2},\Phi_{1}\Phi_{2},\Phi_{1}\Phi_{2}\Phi_{4} & \textrm{ if } q \equiv 0 \mod 2\\
2,\Phi_{2},1/2 \Phi_{1}\Phi_{2},\Phi_{1}\Phi_{2}\Phi_{4} & \textrm{ if } q \equiv 1 \mod 2\\
\end{array}\right.
\\
 & & \textrm{uniform: }\Phi_{1},\Phi_{2},\Phi_{2},\Phi_{1}\Phi_{2}\Phi_{4}\\
A_3 & \begin{minipage}[t]{10em}7\-5\-6\-\\\end{minipage} & \Phi_{1},\Phi_{1},\Phi_{1},\Phi_{1}\Phi_{2}\Phi_{4}\\
D_4(a_1)+A_1 & \begin{minipage}[t]{10em}7\-5\-4\-5\-2\-4\-3\-\\\end{minipage} & \left\{\begin{array}{ll}
\Phi_{1}\Phi_{4},\Phi_{1}\Phi_{2}\Phi_{4} & \textrm{ if } q \equiv 0 \mod 2\\
2,\Phi_{1}\Phi_{4},1/2 \Phi_{1}\Phi_{2}\Phi_{4} & \textrm{ if } q \equiv 1 \mod 2\\
\end{array}\right.
\\
 & & \textrm{uniform: }\Phi_{4},\Phi_{1}\Phi_{2},\Phi_{1}\Phi_{4}\\
A_3+2A_1' & \begin{minipage}[t]{10em}7\-5\-6\-2\-3\-\\\end{minipage} & \Phi_{2},\Phi_{1}\Phi_{2},\Phi_{1}\Phi_{2}\Phi_{4}\\
\end{array}\]
\[
\begin{array}{lll}
\textrm{}&\textrm{}&\textrm{cont. }E_7(q)\\
\hline
D_6+A_1 & \begin{minipage}[t]{10em}7\-6\-5\-4\-5\-2\-4\-5\-3\-4\-5\-2\-4\-3\-1\-\\\end{minipage} & \Phi_{2},\Phi_{2},\Phi_{2}\Phi_{10}\\
A_2+A_1 & \begin{minipage}[t]{10em}7\-6\-4\-\\\end{minipage} & \Phi_{1},\Phi_{1},\Phi_{1},\Phi_{1}\Phi_{2}\Phi_{3}\\
A_2+3A_1 & \begin{minipage}[t]{10em}7\-5\-2\-3\-1\-\\\end{minipage} & \Phi_{2},\Phi_{1}\Phi_{2},\Phi_{1}\Phi_{2}\Phi_{3}\\
A_5+A_2 & \begin{minipage}[t]{10em}7\-6\-5\-6\-4\-5\-6\-2\-4\-3\-4\-5\-6\-2\-4\-5\-3\-1\-3\-4\-5\-2\-4\-3\-1\-\\\end{minipage} & \Phi_{3},\Phi_{2}\Phi_{3}\Phi_{6}\\
D_4+A_1 & \begin{minipage}[t]{10em}7\-5\-4\-2\-3\-\\\end{minipage} & \Phi_{2},\Phi_{1}\Phi_{2},\Phi_{1}\Phi_{2}\Phi_{6}\\
2A_2+A_1 & \begin{minipage}[t]{10em}7\-6\-4\-2\-1\-\\\end{minipage} & \Phi_{1}\Phi_{3},\Phi_{1}\Phi_{2}\Phi_{3}\\
A_5' & \begin{minipage}[t]{10em}7\-5\-2\-6\-4\-\\\end{minipage} & \left\{\begin{array}{ll}
\Phi_{1},\Phi_{1}\Phi_{2}\Phi_{3}\Phi_{6} & \textrm{ if } q \equiv 0 \mod 2\\
2,\Phi_{1},1/2 \Phi_{1}\Phi_{2}\Phi_{3}\Phi_{6} & \textrm{ if } q \equiv 1 \mod 2\\
\end{array}\right.
\\
 & & \textrm{uniform: }\Phi_{1},\Phi_{1}\Phi_{3},\Phi_{2}\Phi_{6}\\
A_5'' & \begin{minipage}[t]{10em}7\-5\-3\-6\-4\-\\\end{minipage} & \Phi_{1},\Phi_{1}\Phi_{2}\Phi_{3}\Phi_{6}\\
E_7(a_1) & \begin{minipage}[t]{10em}7\-6\-5\-4\-5\-2\-4\-3\-1\-\\\end{minipage} & \Phi_{2}\Phi_{14}\\
D_5 & \begin{minipage}[t]{10em}6\-5\-4\-2\-3\-\\\end{minipage} & \Phi_{1},\Phi_{1}\Phi_{2}\Phi_{8}\\
A_7 & \begin{minipage}[t]{10em}7\-6\-7\-5\-6\-4\-5\-2\-4\-5\-3\-4\-5\-2\-4\-3\-1\-\\\end{minipage} & \left\{\begin{array}{ll}
\Phi_{2}\Phi_{4}\Phi_{8} & \textrm{ if } q \equiv 0 \mod 2\\
2,1/2 \Phi_{2}\Phi_{4}\Phi_{8} & \textrm{ if } q \equiv 1 \mod 2\\
\end{array}\right.
\\
 & & \textrm{uniform: }\Phi_{8},\Phi_{2}\Phi_{4}\\
E_7 & \begin{minipage}[t]{10em}7\-6\-5\-4\-2\-3\-1\-\\\end{minipage} & \Phi_{2}\Phi_{18}\\
A_4+A_1 & \begin{minipage}[t]{10em}7\-5\-6\-4\-1\-\\\end{minipage} & \Phi_{1},\Phi_{1}\Phi_{2}\Phi_{5}\\
D_5(a_1) & \begin{minipage}[t]{10em}3\-4\-2\-3\-4\-6\-5\-\\\end{minipage} & \Phi_{1},\Phi_{1}\Phi_{2}\Phi_{4}\Phi_{6}\\
A_3+A_2 & \begin{minipage}[t]{10em}7\-5\-6\-3\-1\-\\\end{minipage} & \Phi_{1},\Phi_{1}\Phi_{2}\Phi_{3}\Phi_{4}\\
E_7(a_2) & \begin{minipage}[t]{10em}7\-6\-7\-5\-6\-4\-5\-2\-4\-3\-1\-\\\end{minipage} & \Phi_{2}\Phi_{6}\Phi_{12}\\
E_7(a_3) & \begin{minipage}[t]{10em}2\-4\-2\-3\-5\-4\-2\-7\-6\-5\-4\-3\-1\-\\\end{minipage} & \Phi_{2}\Phi_{6}\Phi_{10}\\
\end{array}
\]

\newpage
And  here  is the  structure  of  the maximal  tori  in  the simple  group
$(E_7)_{sc}(q)/Z$ in the case $q \equiv 1 \mod 2$ (then $|Z| = 2$).

\[
\begin{array}{lll}
\textrm{class}&\textrm{representative}&\textrm{elementary divisors}\\
\hline
A_0 & \begin{minipage}[t]{10em}-\-\\\end{minipage} & 1/2 \Phi_{1},\Phi_{1},\Phi_{1},\Phi_{1},\Phi_{1},\Phi_{1},\Phi_{1}\\
6A_1 & \begin{minipage}[t]{10em}7\-6\-7\-5\-6\-7\-4\-5\-6\-7\-2\-4\-5\-6\-7\-3\-4\-5\-6\-7\-2\-4\-5\-6\-3\-4\-5\-2\-4\-3\-\\\end{minipage} & 1/2 \Phi_{2},\Phi_{2},\Phi_{2},\Phi_{2},\Phi_{2},\Phi_{1}\Phi_{2}\\
4A_1'' & \begin{minipage}[t]{10em}5\-4\-5\-2\-4\-5\-3\-4\-5\-2\-4\-3\-\\\end{minipage} & \left\{\begin{array}{ll}
2,\Phi_{2},1/4 \Phi_{1}\Phi_{2},\Phi_{1}\Phi_{2},\Phi_{1}\Phi_{2} & \textrm{ if } q \equiv 1 \mod 4\\
\Phi_{2},1/2 \Phi_{1}\Phi_{2},\Phi_{1}\Phi_{2},\Phi_{1}\Phi_{2} & \textrm{ if } q \equiv 3 \mod 4\\
\end{array}\right.
\\
2A_1 & \begin{minipage}[t]{10em}7\-5\-\\\end{minipage} & 1/2 \Phi_{1},\Phi_{1},\Phi_{1},\Phi_{1}\Phi_{2},\Phi_{1}\Phi_{2}\\
4A_1' & \begin{minipage}[t]{10em}7\-5\-2\-3\-\\\end{minipage} & 1/2 \Phi_{2},\Phi_{1}\Phi_{2},\Phi_{1}\Phi_{2},\Phi_{1}\Phi_{2}\\
A_2 & \begin{minipage}[t]{10em}7\-6\-\\\end{minipage} & 1/2 \Phi_{1},\Phi_{1},\Phi_{1},\Phi_{1},\Phi_{1}\Phi_{3}\\
3A_2 & \begin{minipage}[t]{10em}6\-5\-6\-4\-5\-6\-2\-4\-3\-4\-5\-6\-2\-4\-5\-3\-1\-3\-4\-5\-2\-4\-3\-1\-\\\end{minipage} & \Phi_{3},\Phi_{3},1/2 \Phi_{1}\Phi_{3}\\
2A_2 & \begin{minipage}[t]{10em}7\-6\-4\-2\-\\\end{minipage} & 1/2 \Phi_{1},\Phi_{1}\Phi_{3},\Phi_{1}\Phi_{3}\\
D_4(a_1) & \begin{minipage}[t]{10em}5\-4\-5\-2\-4\-3\-\\\end{minipage} & 1/2 \Phi_{1},\Phi_{1}\Phi_{4},\Phi_{1}\Phi_{4}\\
A_3+A_1' & \begin{minipage}[t]{10em}7\-5\-6\-2\-\\\end{minipage} & \left\{\begin{array}{ll}
\Phi_{1},1/2 \Phi_{1}\Phi_{2},\Phi_{1}\Phi_{2}\Phi_{4} & \textrm{ if } q \equiv 1 \mod 4\\
2,\Phi_{1},1/4 \Phi_{1}\Phi_{2},\Phi_{1}\Phi_{2}\Phi_{4} & \textrm{ if } q \equiv 3 \mod 4\\
\end{array}\right.
\\
A_3+3A_1 & \begin{minipage}[t]{10em}7\-5\-4\-5\-2\-4\-5\-3\-4\-5\-2\-4\-3\-1\-\\\end{minipage} & 1/2 \Phi_{2},\Phi_{2},\Phi_{2},\Phi_{1}\Phi_{2}\Phi_{4}\\
D_4(a_1)+2A_1 & \begin{minipage}[t]{10em}7\-6\-7\-5\-6\-4\-5\-2\-4\-5\-3\-4\-5\-2\-4\-3\-\\\end{minipage} & \left\{\begin{array}{ll}
2,1/2 \Phi_{2}\Phi_{4},1/2 \Phi_{1}\Phi_{2}\Phi_{4} & \textrm{ if } q \equiv 1 \mod 4\\
\Phi_{2}\Phi_{4},1/2 \Phi_{1}\Phi_{2}\Phi_{4} & \textrm{ if } q \equiv 3 \mod 4\\
\end{array}\right.
\\
A_3+A_1'' & \begin{minipage}[t]{10em}7\-5\-6\-3\-\\\end{minipage} & 1/2 \Phi_{1},\Phi_{1}\Phi_{2},\Phi_{1}\Phi_{2}\Phi_{4}\\
A_4 & \begin{minipage}[t]{10em}7\-6\-5\-4\-\\\end{minipage} & 1/2 \Phi_{1},\Phi_{1},\Phi_{1}\Phi_{5}\\
D_4+2A_1 & \begin{minipage}[t]{10em}7\-6\-5\-4\-5\-2\-4\-5\-3\-4\-5\-2\-4\-3\-\\\end{minipage} & 1/2 \Phi_{2},\Phi_{2},\Phi_{2},\Phi_{1}\Phi_{2}\Phi_{6}\\
D_4 & \begin{minipage}[t]{10em}5\-4\-2\-3\-\\\end{minipage} & 1/2 \Phi_{1},\Phi_{1}\Phi_{2},\Phi_{1}\Phi_{2}\Phi_{6}\\
E_6(a_2) & \begin{minipage}[t]{10em}1\-2\-3\-1\-5\-4\-6\-5\-4\-2\-3\-4\-\\\end{minipage} & \Phi_{6},1/2 \Phi_{1}\Phi_{3}\Phi_{6}\\
A_2+2A_1 & \begin{minipage}[t]{10em}7\-6\-4\-1\-\\\end{minipage} & 1/2 \Phi_{1},\Phi_{1}\Phi_{2},\Phi_{1}\Phi_{2}\Phi_{3}\\
D_6(a_2) & \begin{minipage}[t]{10em}7\-6\-7\-5\-6\-4\-5\-2\-4\-3\-\\\end{minipage} & 1/2 \Phi_{2}\Phi_{6},\Phi_{1}\Phi_{2}\Phi_{6}\\
A_5+A_1'' & \begin{minipage}[t]{10em}1\-2\-3\-4\-2\-3\-4\-6\-5\-4\-2\-3\-4\-5\-\\\end{minipage} & \left\{\begin{array}{ll}
2,\Phi_{2},1/4 \Phi_{1}\Phi_{2}\Phi_{3}\Phi_{6} & \textrm{ if } q \equiv 1 \mod 4\\
\Phi_{2},1/2 \Phi_{1}\Phi_{2}\Phi_{3}\Phi_{6} & \textrm{ if } q \equiv 3 \mod 4\\
\end{array}\right.
\\
A_5+A_1' & \begin{minipage}[t]{10em}7\-5\-2\-6\-4\-1\-\\\end{minipage} & 1/2 \Phi_{2},\Phi_{1}\Phi_{2}\Phi_{3}\Phi_{6}\\
A_6 & \begin{minipage}[t]{10em}7\-6\-5\-4\-3\-1\-\\\end{minipage} & 1/2 \Phi_{1}\Phi_{7}\\
D_5+A_1 & \begin{minipage}[t]{10em}7\-5\-2\-3\-4\-1\-\\\end{minipage} & 1/2 \Phi_{2},\Phi_{1}\Phi_{2}\Phi_{8}\\
D_6(a_1) & \begin{minipage}[t]{10em}3\-4\-2\-3\-4\-7\-6\-5\-\\\end{minipage} & \left\{\begin{array}{ll}
1/2 \Phi_{1}\Phi_{4}\Phi_{8} & \textrm{ if } q \equiv 1 \mod 4\\
2,1/4 \Phi_{1}\Phi_{4}\Phi_{8} & \textrm{ if } q \equiv 3 \mod 4\\
\end{array}\right.
\\
\end{array}\]
\[
\begin{array}{lll}
\textrm{}&\textrm{}&\textrm{cont. }E_7(q)/Z\\
\hline
E_6(a_1) & \begin{minipage}[t]{10em}6\-5\-4\-5\-2\-4\-3\-1\-\\\end{minipage} & 1/2 \Phi_{1}\Phi_{9}\\
D_6 & \begin{minipage}[t]{10em}7\-6\-5\-4\-2\-3\-\\\end{minipage} & 1/2 \Phi_{2},\Phi_{1}\Phi_{2}\Phi_{10}\\
A_3+A_2+A_1 & \begin{minipage}[t]{10em}7\-5\-6\-2\-3\-1\-\\\end{minipage} & 1/2 \Phi_{2},\Phi_{1}\Phi_{2}\Phi_{3}\Phi_{4}\\
D_5(a_1)+A_1 & \begin{minipage}[t]{10em}7\-5\-4\-5\-2\-4\-3\-1\-\\\end{minipage} & 1/2 \Phi_{2},\Phi_{1}\Phi_{2}\Phi_{4}\Phi_{6}\\
E_6 & \begin{minipage}[t]{10em}6\-4\-1\-5\-3\-2\-\\\end{minipage} & 1/2 \Phi_{1}\Phi_{3}\Phi_{12}\\
A_4+A_2 & \begin{minipage}[t]{10em}7\-6\-4\-2\-3\-1\-\\\end{minipage} & 1/2 \Phi_{1}\Phi_{3}\Phi_{5}\\
7A_1 & \begin{minipage}[t]{10em}7\-6\-7\-5\-6\-7\-4\-5\-6\-7\-2\-4\-5\-6\-7\-3\-4\-5\-6\-7\-2\-4\-5\-6\-3\-4\-5\-2\-4\-3\-1\-3\-4\-5\-6\-7\-2\-4\-5\-6\-3\-4\-5\-2\-4\-3\-1\-3\-4\-5\-6\-7\-2\-4\-5\-6\-3\-4\-5\-2\-4\-3\-1\-\\\end{minipage} & 1/2 \Phi_{2},\Phi_{2},\Phi_{2},\Phi_{2},\Phi_{2},\Phi_{2},\Phi_{2}\\
A_1 & \begin{minipage}[t]{10em}7\-\\\end{minipage} & 1/2 \Phi_{1},\Phi_{1},\Phi_{1},\Phi_{1},\Phi_{1},\Phi_{1}\Phi_{2}\\
3A_1' & \begin{minipage}[t]{10em}7\-5\-2\-\\\end{minipage} & \left\{\begin{array}{ll}
\Phi_{1},1/2 \Phi_{1}\Phi_{2},\Phi_{1}\Phi_{2},\Phi_{1}\Phi_{2} & \textrm{ if } q \equiv 1 \mod 4\\
2,\Phi_{1},1/4 \Phi_{1}\Phi_{2},\Phi_{1}\Phi_{2},\Phi_{1}\Phi_{2} & \textrm{ if } q \equiv 3 \mod 4\\
\end{array}\right.
\\
5A_1 & \begin{minipage}[t]{10em}7\-5\-4\-5\-2\-4\-5\-3\-4\-5\-2\-4\-3\-\\\end{minipage} & 1/2 \Phi_{2},\Phi_{2},\Phi_{2},\Phi_{1}\Phi_{2},\Phi_{1}\Phi_{2}\\
3A_1'' & \begin{minipage}[t]{10em}7\-5\-3\-\\\end{minipage} & 1/2 \Phi_{1},\Phi_{1}\Phi_{2},\Phi_{1}\Phi_{2},\Phi_{1}\Phi_{2}\\
D_4+3A_1 & \begin{minipage}[t]{10em}7\-6\-7\-5\-6\-7\-4\-5\-6\-7\-2\-4\-5\-6\-7\-3\-4\-5\-6\-7\-2\-4\-5\-6\-3\-4\-5\-2\-4\-3\-1\-\\\end{minipage} & 1/2 \Phi_{2},\Phi_{2},\Phi_{2},\Phi_{2},\Phi_{2}\Phi_{6}\\
E_7(a_4) & \begin{minipage}[t]{10em}7\-6\-7\-5\-6\-7\-4\-5\-6\-2\-4\-3\-4\-5\-2\-4\-3\-1\-3\-4\-2\-\\\end{minipage} & \Phi_{6},\Phi_{6},1/2 \Phi_{2}\Phi_{6}\\
D_6(a_2)+A_1 & \begin{minipage}[t]{10em}7\-6\-7\-5\-6\-4\-5\-6\-2\-4\-5\-6\-3\-4\-5\-6\-1\-3\-4\-5\-2\-4\-3\-\\\end{minipage} & 1/2 \Phi_{2},\Phi_{2}\Phi_{6},\Phi_{2}\Phi_{6}\\
2A_3+A_1 & \begin{minipage}[t]{10em}7\-6\-7\-5\-6\-7\-4\-5\-6\-7\-2\-4\-5\-6\-7\-3\-4\-5\-6\-7\-2\-4\-5\-6\-3\-1\-3\-4\-5\-2\-4\-3\-1\-\\\end{minipage} & 1/2 \Phi_{2},\Phi_{2}\Phi_{4},\Phi_{2}\Phi_{4}\\
A_3+2A_1'' & \begin{minipage}[t]{10em}6\-5\-4\-5\-2\-4\-5\-3\-4\-5\-2\-4\-3\-\\\end{minipage} & \left\{\begin{array}{ll}
2,\Phi_{2},1/4 \Phi_{1}\Phi_{2},\Phi_{1}\Phi_{2}\Phi_{4} & \textrm{ if } q \equiv 1 \mod 4\\
\Phi_{2},1/2 \Phi_{1}\Phi_{2},\Phi_{1}\Phi_{2}\Phi_{4} & \textrm{ if } q \equiv 3 \mod 4\\
\end{array}\right.
\\
A_3 & \begin{minipage}[t]{10em}7\-5\-6\-\\\end{minipage} & 1/2 \Phi_{1},\Phi_{1},\Phi_{1},\Phi_{1}\Phi_{2}\Phi_{4}\\
D_4(a_1)+A_1 & \begin{minipage}[t]{10em}7\-5\-4\-5\-2\-4\-3\-\\\end{minipage} & \left\{\begin{array}{ll}
\Phi_{1}\Phi_{4},1/2 \Phi_{1}\Phi_{2}\Phi_{4} & \textrm{ if } q \equiv 1 \mod 4\\
2,1/2 \Phi_{1}\Phi_{4},1/2 \Phi_{1}\Phi_{2}\Phi_{4} & \textrm{ if } q \equiv 3 \mod 4\\
\end{array}\right.
\\
A_3+2A_1' & \begin{minipage}[t]{10em}7\-5\-6\-2\-3\-\\\end{minipage} & 1/2 \Phi_{2},\Phi_{1}\Phi_{2},\Phi_{1}\Phi_{2}\Phi_{4}\\
D_6+A_1 & \begin{minipage}[t]{10em}7\-6\-5\-4\-5\-2\-4\-5\-3\-4\-5\-2\-4\-3\-1\-\\\end{minipage} & 1/2 \Phi_{2},\Phi_{2},\Phi_{2}\Phi_{10}\\
A_2+A_1 & \begin{minipage}[t]{10em}7\-6\-4\-\\\end{minipage} & 1/2 \Phi_{1},\Phi_{1},\Phi_{1},\Phi_{1}\Phi_{2}\Phi_{3}\\
A_2+3A_1 & \begin{minipage}[t]{10em}7\-5\-2\-3\-1\-\\\end{minipage} & 1/2 \Phi_{2},\Phi_{1}\Phi_{2},\Phi_{1}\Phi_{2}\Phi_{3}\\
A_5+A_2 & \begin{minipage}[t]{10em}7\-6\-5\-6\-4\-5\-6\-2\-4\-3\-4\-5\-6\-2\-4\-5\-3\-1\-3\-4\-5\-2\-4\-3\-1\-\\\end{minipage} & \Phi_{3},1/2 \Phi_{2}\Phi_{3}\Phi_{6}\\
\end{array}\]
\[
\begin{array}{lll}
\textrm{}&\textrm{}&\textrm{cont. }E_7(q)/Z\\
\hline
D_4+A_1 & \begin{minipage}[t]{10em}7\-5\-4\-2\-3\-\\\end{minipage} & 1/2 \Phi_{2},\Phi_{1}\Phi_{2},\Phi_{1}\Phi_{2}\Phi_{6}\\
2A_2+A_1 & \begin{minipage}[t]{10em}7\-6\-4\-2\-1\-\\\end{minipage} & 1/2 \Phi_{1}\Phi_{3},\Phi_{1}\Phi_{2}\Phi_{3}\\
A_5' & \begin{minipage}[t]{10em}7\-5\-2\-6\-4\-\\\end{minipage} & \left\{\begin{array}{ll}
\Phi_{1},1/2 \Phi_{1}\Phi_{2}\Phi_{3}\Phi_{6} & \textrm{ if } q \equiv 1 \mod 4\\
2,\Phi_{1},1/4 \Phi_{1}\Phi_{2}\Phi_{3}\Phi_{6} & \textrm{ if } q \equiv 3 \mod 4\\
\end{array}\right.
\\
A_5'' & \begin{minipage}[t]{10em}7\-5\-3\-6\-4\-\\\end{minipage} & 1/2 \Phi_{1},\Phi_{1}\Phi_{2}\Phi_{3}\Phi_{6}\\
E_7(a_1) & \begin{minipage}[t]{10em}7\-6\-5\-4\-5\-2\-4\-3\-1\-\\\end{minipage} & 1/2 \Phi_{2}\Phi_{14}\\
D_5 & \begin{minipage}[t]{10em}6\-5\-4\-2\-3\-\\\end{minipage} & 1/2 \Phi_{1},\Phi_{1}\Phi_{2}\Phi_{8}\\
A_7 & \begin{minipage}[t]{10em}7\-6\-7\-5\-6\-4\-5\-2\-4\-5\-3\-4\-5\-2\-4\-3\-1\-\\\end{minipage} & \left\{\begin{array}{ll}
2,1/4 \Phi_{2}\Phi_{4}\Phi_{8} & \textrm{ if } q \equiv 1 \mod 4\\
1/2 \Phi_{2}\Phi_{4}\Phi_{8} & \textrm{ if } q \equiv 3 \mod 4\\
\end{array}\right.
\\
E_7 & \begin{minipage}[t]{10em}7\-6\-5\-4\-2\-3\-1\-\\\end{minipage} & 1/2 \Phi_{2}\Phi_{18}\\
A_4+A_1 & \begin{minipage}[t]{10em}7\-5\-6\-4\-1\-\\\end{minipage} & 1/2 \Phi_{1},\Phi_{1}\Phi_{2}\Phi_{5}\\
D_5(a_1) & \begin{minipage}[t]{10em}3\-4\-2\-3\-4\-6\-5\-\\\end{minipage} & 1/2 \Phi_{1},\Phi_{1}\Phi_{2}\Phi_{4}\Phi_{6}\\
A_3+A_2 & \begin{minipage}[t]{10em}7\-5\-6\-3\-1\-\\\end{minipage} & 1/2 \Phi_{1},\Phi_{1}\Phi_{2}\Phi_{3}\Phi_{4}\\
E_7(a_2) & \begin{minipage}[t]{10em}7\-6\-7\-5\-6\-4\-5\-2\-4\-3\-1\-\\\end{minipage} & 1/2 \Phi_{2}\Phi_{6}\Phi_{12}\\
E_7(a_3) & \begin{minipage}[t]{10em}2\-4\-2\-3\-5\-4\-2\-7\-6\-5\-4\-3\-1\-\\\end{minipage} & 1/2 \Phi_{2}\Phi_{6}\Phi_{10}\\
\end{array}
\]

\newpage
\subsection{$E_8(q)$}

\[
\begin{array}{lll}
\textrm{class}&\textrm{representative}&\textrm{elementary divisors}\\
\hline
A_0 & \begin{minipage}[t]{8em}-\-\\\end{minipage} & \Phi_{1},\Phi_{1},\Phi_{1},\Phi_{1},\Phi_{1},\Phi_{1},\Phi_{1},\Phi_{1}\\
8A_1 & \begin{minipage}[t]{8em}8\-7\-8\-6\-7\-8\-5\-6\-7\-8\-4\-5\-6\-7\-8\-2\-4\-5\-6\-7\-8\-3\-4\-5\-6\-7\-8\-2\-4\-5\-6\-7\-3\-4\-5\-6\-2\-4\-5\-3\-4\-2\-1\-3\-4\-5\-6\-7\-8\-2\-4\-5\-6\-7\-3\-4\-5\-6\-2\-4\-5\-3\-4\-2\-1\-3\-4\-5\-6\-7\-8\-2\-4\-5\-6\-7\-3\-4\-5\-6\-2\-4\-5\-3\-4\-2\-1\-3\-4\-5\-6\-7\-8\-2\-4\-5\-6\-3\-4\-5\-2\-4\-3\-1\-3\-4\-5\-6\-7\-2\-4\-5\-6\-3\-4\-5\-2\-4\-3\-1\-\\\end{minipage} & \Phi_{2},\Phi_{2},\Phi_{2},\Phi_{2},\Phi_{2},\Phi_{2},\Phi_{2},\Phi_{2}\\
4A_1' & \begin{minipage}[t]{8em}5\-4\-5\-2\-4\-5\-3\-4\-5\-2\-4\-3\-\\\end{minipage} & \left\{\begin{array}{ll}
\Phi_{1}\Phi_{2},\Phi_{1}\Phi_{2},\Phi_{1}\Phi_{2},\Phi_{1}\Phi_{2} & \textrm{ if } q \equiv 0 \mod 2\\
2,2,1/2 \Phi_{1}\Phi_{2},1/2 \Phi_{1}\Phi_{2},\Phi_{1}\Phi_{2},\Phi_{1}\Phi_{2} & \textrm{ if } q \equiv 1 \mod 2\\
\end{array}\right.
\\
 & & \textrm{uniform: }\Phi_{1},\Phi_{2},\Phi_{1},\Phi_{2},\Phi_{1}\Phi_{2},\Phi_{1}\Phi_{2}\\
2A_1 & \begin{minipage}[t]{8em}6\-1\-\\\end{minipage} & \Phi_{1},\Phi_{1},\Phi_{1},\Phi_{1},\Phi_{1}\Phi_{2},\Phi_{1}\Phi_{2}\\
6A_1 & \begin{minipage}[t]{8em}7\-6\-7\-5\-6\-7\-4\-5\-6\-7\-2\-4\-5\-6\-7\-3\-4\-5\-6\-7\-2\-4\-5\-6\-3\-4\-5\-2\-4\-3\-\\\end{minipage} & \Phi_{2},\Phi_{2},\Phi_{2},\Phi_{2},\Phi_{1}\Phi_{2},\Phi_{1}\Phi_{2}\\
2D_4(a_1) & \begin{minipage}[t]{8em}8\-7\-6\-7\-5\-6\-4\-5\-6\-7\-2\-4\-5\-6\-3\-4\-5\-6\-7\-8\-2\-4\-5\-6\-3\-4\-5\-1\-3\-4\-5\-6\-7\-8\-2\-4\-5\-6\-7\-3\-4\-5\-6\-2\-4\-5\-3\-4\-2\-1\-3\-4\-5\-6\-7\-2\-4\-5\-3\-1\-\\\end{minipage} & \Phi_{4},\Phi_{4},\Phi_{4},\Phi_{4}\\
4A_1'' & \begin{minipage}[t]{8em}7\-5\-2\-3\-\\\end{minipage} & \Phi_{1}\Phi_{2},\Phi_{1}\Phi_{2},\Phi_{1}\Phi_{2},\Phi_{1}\Phi_{2}\\
A_2 & \begin{minipage}[t]{8em}6\-7\-\\\end{minipage} & \Phi_{1},\Phi_{1},\Phi_{1},\Phi_{1},\Phi_{1},\Phi_{1}\Phi_{3}\\
D_4+4A_1 & \begin{minipage}[t]{8em}1\-2\-3\-1\-4\-2\-3\-1\-4\-3\-5\-4\-2\-3\-1\-4\-3\-5\-4\-2\-6\-5\-4\-2\-3\-1\-4\-3\-5\-4\-2\-6\-5\-4\-3\-1\-7\-6\-5\-4\-2\-3\-1\-4\-3\-5\-4\-2\-6\-5\-4\-3\-1\-7\-6\-5\-4\-2\-3\-4\-5\-6\-7\-8\-\\\end{minipage} & \Phi_{2},\Phi_{2},\Phi_{2},\Phi_{2},\Phi_{2},\Phi_{2}\Phi_{6}\\
4A_2 & \begin{minipage}[t]{8em}8\-7\-8\-6\-7\-5\-6\-7\-4\-2\-4\-5\-3\-4\-5\-6\-7\-8\-2\-4\-5\-6\-3\-4\-5\-2\-1\-3\-4\-5\-6\-7\-8\-2\-4\-5\-6\-7\-3\-4\-5\-6\-2\-4\-5\-3\-4\-1\-3\-4\-5\-6\-7\-8\-2\-4\-5\-6\-7\-3\-4\-5\-6\-2\-4\-5\-3\-4\-2\-1\-3\-4\-5\-6\-2\-4\-5\-3\-4\-1\-\\\end{minipage} & \Phi_{3},\Phi_{3},\Phi_{3},\Phi_{3}\\
E_8(a_8) & \begin{minipage}[t]{8em}5\-6\-4\-3\-4\-5\-6\-7\-8\-2\-4\-5\-6\-3\-4\-2\-1\-3\-4\-5\-6\-7\-8\-2\-4\-5\-6\-7\-3\-4\-5\-2\-1\-3\-4\-5\-6\-7\-8\-2\-\\\end{minipage} & \Phi_{6},\Phi_{6},\Phi_{6},\Phi_{6}\\
3A_2 & \begin{minipage}[t]{8em}6\-5\-2\-4\-5\-6\-3\-4\-5\-2\-4\-3\-1\-3\-4\-5\-6\-2\-4\-5\-3\-4\-1\-3\-\\\end{minipage} & \left\{\begin{array}{ll}
\Phi_{3},\Phi_{1}\Phi_{3},\Phi_{1}\Phi_{3} & \textrm{ if } q \equiv 0,2 \mod 3\\
3,\Phi_{3},1/3 \Phi_{1}\Phi_{3},\Phi_{1}\Phi_{3} & \textrm{ if } q \equiv 1 \mod 3\\
\end{array}\right.
\\
 & & \textrm{uniform: }\Phi_{1},\Phi_{3},\Phi_{3},\Phi_{1}\Phi_{3}\\
\end{array}\]
\[
\begin{array}{lll}
\textrm{}&\textrm{}&\textrm{cont. }E_8(q)\\
\hline
E_7(a_4)+A_1 & \begin{minipage}[t]{8em}2\-3\-4\-2\-3\-4\-5\-4\-2\-3\-1\-4\-5\-6\-8\-7\-6\-5\-4\-2\-3\-1\-4\-3\-5\-4\-2\-6\-5\-4\-3\-1\-7\-6\-5\-4\-2\-3\-4\-5\-6\-7\-\\\end{minipage} & \left\{\begin{array}{ll}
\Phi_{6},\Phi_{2}\Phi_{6},\Phi_{2}\Phi_{6} & \textrm{ if } q \equiv 0,1 \mod 3\\
3,\Phi_{6},1/3 \Phi_{2}\Phi_{6},\Phi_{2}\Phi_{6} & \textrm{ if } q \equiv 2 \mod 3\\
\end{array}\right.
\\
 & & \textrm{uniform: }\Phi_{2},\Phi_{6},\Phi_{6},\Phi_{2}\Phi_{6}\\
2A_2 & \begin{minipage}[t]{8em}7\-8\-2\-4\-\\\end{minipage} & \Phi_{1},\Phi_{1},\Phi_{1}\Phi_{3},\Phi_{1}\Phi_{3}\\
2D_4 & \begin{minipage}[t]{8em}2\-4\-2\-5\-4\-2\-6\-5\-4\-2\-3\-4\-5\-6\-7\-6\-5\-4\-2\-3\-4\-5\-6\-7\-8\-7\-6\-5\-4\-2\-3\-1\-4\-3\-5\-4\-2\-6\-5\-4\-3\-1\-7\-8\-\\\end{minipage} & \Phi_{2},\Phi_{2},\Phi_{2}\Phi_{6},\Phi_{2}\Phi_{6}\\
D_4(a_1) & \begin{minipage}[t]{8em}4\-2\-4\-3\-4\-5\-\\\end{minipage} & \Phi_{1},\Phi_{1},\Phi_{1}\Phi_{4},\Phi_{1}\Phi_{4}\\
2A_3+2A_1 & \begin{minipage}[t]{8em}1\-2\-3\-1\-4\-2\-3\-1\-4\-3\-5\-4\-2\-3\-1\-4\-3\-5\-4\-2\-6\-5\-4\-2\-3\-1\-4\-3\-5\-4\-2\-6\-5\-4\-3\-1\-7\-6\-5\-4\-2\-3\-1\-4\-3\-5\-4\-2\-6\-5\-4\-3\-1\-7\-8\-7\-6\-5\-4\-2\-3\-4\-5\-6\-7\-8\-\\\end{minipage} & \Phi_{2},\Phi_{2},\Phi_{2}\Phi_{4},\Phi_{2}\Phi_{4}\\
2A_3' & \begin{minipage}[t]{8em}2\-3\-4\-2\-3\-4\-6\-5\-7\-6\-5\-4\-2\-3\-4\-5\-\\\end{minipage} & \left\{\begin{array}{ll}
\Phi_{1}\Phi_{2}\Phi_{4},\Phi_{1}\Phi_{2}\Phi_{4} & \textrm{ if } q \equiv 0 \mod 2\\
2,2,1/2 \Phi_{1}\Phi_{2}\Phi_{4},1/2 \Phi_{1}\Phi_{2}\Phi_{4} & \textrm{ if } q \equiv 1 \mod 2\\
\end{array}\right.
\\
 & & \textrm{uniform: }\Phi_{4},\Phi_{4},\Phi_{1}\Phi_{2},\Phi_{1}\Phi_{2}\\
A_3+A_1 & \begin{minipage}[t]{8em}3\-7\-6\-8\-\\\end{minipage} & \Phi_{1},\Phi_{1},\Phi_{1}\Phi_{2},\Phi_{1}\Phi_{2}\Phi_{4}\\
A_3+3A_1 & \begin{minipage}[t]{8em}1\-2\-3\-4\-2\-3\-4\-5\-4\-2\-3\-4\-5\-7\-\\\end{minipage} & \Phi_{2},\Phi_{2},\Phi_{1}\Phi_{2},\Phi_{1}\Phi_{2}\Phi_{4}\\
D_8(a_3) & \begin{minipage}[t]{8em}4\-5\-6\-7\-3\-4\-5\-2\-4\-3\-1\-3\-4\-5\-6\-7\-8\-2\-4\-5\-3\-4\-2\-1\-3\-4\-5\-6\-7\-8\-\\\end{minipage} & \Phi_{8},\Phi_{8}\\
2A_3'' & \begin{minipage}[t]{8em}7\-6\-8\-1\-4\-3\-\\\end{minipage} & \Phi_{1}\Phi_{2}\Phi_{4},\Phi_{1}\Phi_{2}\Phi_{4}\\
A_4 & \begin{minipage}[t]{8em}5\-6\-3\-4\-\\\end{minipage} & \Phi_{1},\Phi_{1},\Phi_{1},\Phi_{1}\Phi_{5}\\
D_6+2A_1 & \begin{minipage}[t]{8em}3\-4\-3\-5\-4\-3\-6\-5\-4\-3\-8\-7\-6\-5\-4\-2\-3\-1\-4\-3\-5\-4\-2\-6\-5\-4\-3\-7\-6\-5\-4\-2\-\\\end{minipage} & \Phi_{2},\Phi_{2},\Phi_{2},\Phi_{2}\Phi_{10}\\
2A_4 & \begin{minipage}[t]{8em}7\-5\-6\-2\-4\-5\-3\-4\-5\-6\-7\-8\-2\-4\-5\-6\-7\-1\-3\-4\-5\-6\-7\-8\-2\-4\-5\-6\-7\-3\-4\-5\-2\-4\-3\-1\-3\-4\-5\-6\-7\-8\-2\-4\-5\-6\-3\-4\-\\\end{minipage} & \Phi_{5},\Phi_{5}\\
E_8(a_6) & \begin{minipage}[t]{8em}8\-7\-6\-7\-5\-4\-5\-6\-2\-4\-3\-4\-5\-6\-7\-1\-3\-4\-5\-6\-2\-4\-5\-3\-\\\end{minipage} & \Phi_{10},\Phi_{10}\\
A_2+4A_1 & \begin{minipage}[t]{8em}8\-7\-5\-4\-5\-2\-4\-5\-3\-4\-5\-2\-4\-3\-\\\end{minipage} & \Phi_{2},\Phi_{2},\Phi_{1}\Phi_{2},\Phi_{1}\Phi_{2}\Phi_{3}\\
D_4 & \begin{minipage}[t]{8em}5\-4\-2\-3\-\\\end{minipage} & \Phi_{1},\Phi_{1},\Phi_{1}\Phi_{2},\Phi_{1}\Phi_{2}\Phi_{6}\\
E_6(a_2)+A_2 & \begin{minipage}[t]{8em}1\-2\-3\-1\-4\-2\-3\-1\-4\-5\-4\-2\-3\-1\-4\-3\-6\-5\-4\-2\-3\-1\-4\-3\-5\-4\-2\-6\-5\-4\-3\-1\-7\-6\-5\-4\-2\-3\-4\-5\-6\-7\-8\-7\-\\\end{minipage} & \Phi_{3}\Phi_{6},\Phi_{3}\Phi_{6}\\
\end{array}\]
\[
\begin{array}{lll}
\textrm{}&\textrm{}&\textrm{cont. }E_8(q)\\
\hline
A_2+2A_1 & \begin{minipage}[t]{8em}7\-4\-5\-1\-\\\end{minipage} & \Phi_{1},\Phi_{1},\Phi_{1}\Phi_{2},\Phi_{1}\Phi_{2}\Phi_{3}\\
D_4+2A_1 & \begin{minipage}[t]{8em}5\-4\-5\-2\-4\-5\-3\-4\-5\-6\-7\-2\-4\-3\-\\\end{minipage} & \Phi_{2},\Phi_{2},\Phi_{1}\Phi_{2},\Phi_{1}\Phi_{2}\Phi_{6}\\
E_8(a_3) & \begin{minipage}[t]{8em}8\-7\-6\-7\-5\-2\-4\-3\-4\-5\-2\-4\-1\-3\-4\-5\-6\-2\-4\-5\-\\\end{minipage} & \Phi_{12},\Phi_{12}\\
E_6(a_2) & \begin{minipage}[t]{8em}3\-4\-2\-3\-5\-4\-2\-3\-1\-4\-5\-6\-\\\end{minipage} & \Phi_{1}\Phi_{6},\Phi_{1}\Phi_{3}\Phi_{6}\\
A_5+A_2+A_1 & \begin{minipage}[t]{8em}4\-3\-5\-4\-2\-3\-4\-5\-6\-5\-4\-2\-3\-4\-5\-6\-7\-6\-5\-4\-2\-3\-4\-5\-6\-7\-8\-7\-6\-5\-4\-2\-3\-1\-4\-3\-5\-4\-2\-6\-5\-4\-3\-1\-7\-8\-\\\end{minipage} & \Phi_{2}\Phi_{3},\Phi_{2}\Phi_{3}\Phi_{6}\\
A_5+A_1' & \begin{minipage}[t]{8em}2\-3\-4\-2\-3\-4\-5\-4\-2\-3\-1\-4\-5\-6\-\\\end{minipage} & \left\{\begin{array}{ll}
\Phi_{1}\Phi_{2},\Phi_{1}\Phi_{2}\Phi_{3}\Phi_{6} & \textrm{ if } q \equiv 0 \mod 2\\
2,2,1/2 \Phi_{1}\Phi_{2},1/2 \Phi_{1}\Phi_{2}\Phi_{3}\Phi_{6} & \textrm{ if } q \equiv 1 \mod 2\\
\end{array}\right.
\\
 & & \textrm{uniform: }\Phi_{1},\Phi_{2},\Phi_{1}\Phi_{3},\Phi_{2}\Phi_{6}\\
D_4+A_2 & \begin{minipage}[t]{8em}8\-7\-5\-2\-4\-3\-\\\end{minipage} & \Phi_{1}\Phi_{2},\Phi_{1}\Phi_{2}\Phi_{3}\Phi_{6}\\
2A_2+2A_1 & \begin{minipage}[t]{8em}7\-8\-5\-2\-3\-1\-\\\end{minipage} & \Phi_{1}\Phi_{2}\Phi_{3},\Phi_{1}\Phi_{2}\Phi_{3}\\
D_6(a_2) & \begin{minipage}[t]{8em}6\-5\-6\-7\-4\-2\-4\-5\-3\-4\-\\\end{minipage} & \Phi_{1}\Phi_{2}\Phi_{6},\Phi_{1}\Phi_{2}\Phi_{6}\\
D_8(a_1) & \begin{minipage}[t]{8em}5\-4\-6\-5\-4\-2\-3\-7\-6\-5\-4\-2\-3\-8\-7\-6\-5\-4\-2\-3\-1\-4\-\\\end{minipage} & \Phi_{4},\Phi_{4}\Phi_{12}\\
A_5+A_1'' & \begin{minipage}[t]{8em}6\-4\-5\-2\-7\-1\-\\\end{minipage} & \Phi_{1}\Phi_{2},\Phi_{1}\Phi_{2}\Phi_{3}\Phi_{6}\\
A_6 & \begin{minipage}[t]{8em}8\-5\-6\-7\-2\-4\-\\\end{minipage} & \Phi_{1},\Phi_{1}\Phi_{7}\\
D_8 & \begin{minipage}[t]{8em}2\-7\-6\-5\-4\-2\-3\-4\-8\-7\-6\-5\-4\-2\-3\-1\-4\-3\-\\\end{minipage} & \Phi_{2},\Phi_{2}\Phi_{14}\\
D_5+A_1 & \begin{minipage}[t]{8em}8\-5\-6\-2\-3\-4\-\\\end{minipage} & \Phi_{1}\Phi_{2},\Phi_{1}\Phi_{2}\Phi_{8}\\
D_6(a_1) & \begin{minipage}[t]{8em}2\-4\-2\-3\-4\-5\-6\-7\-\\\end{minipage} & \left\{\begin{array}{ll}
\Phi_{1},\Phi_{1}\Phi_{4}\Phi_{8} & \textrm{ if } q \equiv 0 \mod 2\\
2 \Phi_{1},1/2 \Phi_{1}\Phi_{4}\Phi_{8} & \textrm{ if } q \equiv 1 \mod 2\\
\end{array}\right.
\\
 & & \textrm{uniform: }\Phi_{1}\Phi_{4},\Phi_{1}\Phi_{8}\\
A_7+A_1 & \begin{minipage}[t]{8em}2\-4\-2\-5\-4\-2\-3\-6\-5\-4\-2\-3\-4\-5\-6\-7\-6\-5\-4\-2\-3\-1\-4\-3\-5\-4\-2\-6\-5\-4\-3\-1\-7\-8\-\\\end{minipage} & \left\{\begin{array}{ll}
\Phi_{2},\Phi_{2}\Phi_{4}\Phi_{8} & \textrm{ if } q \equiv 0 \mod 2\\
2 \Phi_{2},1/2 \Phi_{2}\Phi_{4}\Phi_{8} & \textrm{ if } q \equiv 1 \mod 2\\
\end{array}\right.
\\
 & & \textrm{uniform: }\Phi_{2}\Phi_{4},\Phi_{2}\Phi_{8}\\
E_6(a_1) & \begin{minipage}[t]{8em}5\-3\-4\-5\-6\-2\-4\-1\-\\\end{minipage} & \Phi_{1},\Phi_{1}\Phi_{9}\\
E_7+A_1 & \begin{minipage}[t]{8em}3\-4\-3\-5\-4\-2\-3\-1\-4\-3\-5\-4\-2\-8\-7\-6\-\\\end{minipage} & \Phi_{2},\Phi_{2}\Phi_{18}\\
A_8 & \begin{minipage}[t]{8em}1\-3\-4\-3\-1\-5\-4\-2\-3\-4\-5\-8\-7\-6\-5\-4\-2\-3\-1\-4\-3\-5\-4\-2\-6\-5\-7\-6\-\\\end{minipage} & \left\{\begin{array}{ll}
\Phi_{3}\Phi_{9} & \textrm{ if } q \equiv 0,2 \mod 3\\
3,1/3 \Phi_{3}\Phi_{9} & \textrm{ if } q \equiv 1 \mod 3\\
\end{array}\right.
\\
 & & \textrm{uniform: }\Phi_{3},\Phi_{9}\\
\end{array}\]
\[
\begin{array}{lll}
\textrm{}&\textrm{}&\textrm{cont. }E_8(q)\\
\hline
E_8(a_4) & \begin{minipage}[t]{8em}2\-3\-4\-3\-1\-5\-4\-2\-3\-4\-5\-6\-7\-8\-\\\end{minipage} & \left\{\begin{array}{ll}
\Phi_{6}\Phi_{18} & \textrm{ if } q \equiv 0,1 \mod 3\\
3,1/3 \Phi_{6}\Phi_{18} & \textrm{ if } q \equiv 2 \mod 3\\
\end{array}\right.
\\
 & & \textrm{uniform: }\Phi_{6},\Phi_{18}\\
A_4+2A_1 & \begin{minipage}[t]{8em}6\-7\-8\-5\-2\-1\-\\\end{minipage} & \Phi_{1}\Phi_{2},\Phi_{1}\Phi_{2}\Phi_{5}\\
D_6 & \begin{minipage}[t]{8em}7\-5\-6\-2\-3\-4\-\\\end{minipage} & \Phi_{1}\Phi_{2},\Phi_{1}\Phi_{2}\Phi_{10}\\
E_8(a_2) & \begin{minipage}[t]{8em}6\-7\-8\-5\-6\-7\-4\-5\-3\-4\-2\-1\-\\\end{minipage} & \Phi_{20}\\
D_4(a_1)+A_2 & \begin{minipage}[t]{8em}7\-8\-4\-2\-3\-4\-5\-2\-\\\end{minipage} & \Phi_{1}\Phi_{4},\Phi_{1}\Phi_{3}\Phi_{4}\\
D_5(a_1)+A_3 & \begin{minipage}[t]{8em}1\-3\-1\-4\-5\-4\-3\-6\-5\-4\-2\-3\-1\-4\-3\-5\-4\-2\-7\-6\-5\-4\-2\-3\-1\-4\-3\-5\-4\-2\-6\-5\-4\-3\-1\-7\-6\-5\-4\-2\-3\-4\-5\-6\-7\-8\-\\\end{minipage} & \Phi_{2}\Phi_{4},\Phi_{2}\Phi_{4}\Phi_{6}\\
E_6+A_2 & \begin{minipage}[t]{8em}1\-2\-3\-4\-2\-3\-1\-4\-3\-5\-4\-2\-3\-1\-4\-5\-6\-5\-4\-2\-3\-4\-5\-6\-8\-7\-\\\end{minipage} & \Phi_{3},\Phi_{3}\Phi_{12}\\
E_8(a_7) & \begin{minipage}[t]{8em}1\-2\-3\-1\-4\-2\-5\-4\-2\-3\-1\-4\-6\-5\-4\-2\-3\-4\-5\-6\-7\-8\-\\\end{minipage} & \Phi_{6},\Phi_{6}\Phi_{12}\\
E_6 & \begin{minipage}[t]{8em}6\-4\-1\-5\-3\-2\-\\\end{minipage} & \Phi_{1},\Phi_{1}\Phi_{3}\Phi_{12}\\
E_7(a_2)+A_1 & \begin{minipage}[t]{8em}1\-2\-3\-1\-4\-2\-3\-1\-4\-3\-5\-4\-3\-1\-6\-5\-4\-2\-3\-4\-5\-6\-7\-8\-\\\end{minipage} & \Phi_{2},\Phi_{2}\Phi_{6}\Phi_{12}\\
A_3+A_2+A_1 & \begin{minipage}[t]{8em}5\-7\-6\-2\-1\-3\-\\\end{minipage} & \Phi_{1}\Phi_{2},\Phi_{1}\Phi_{2}\Phi_{3}\Phi_{4}\\
D_5(a_1)+A_1 & \begin{minipage}[t]{8em}8\-1\-2\-3\-4\-2\-5\-4\-\\\end{minipage} & \Phi_{1}\Phi_{2},\Phi_{1}\Phi_{2}\Phi_{4}\Phi_{6}\\
E_8(a_1) & \begin{minipage}[t]{8em}8\-7\-6\-2\-4\-5\-3\-4\-1\-3\-\\\end{minipage} & \Phi_{24}\\
A_4+A_2 & \begin{minipage}[t]{8em}8\-7\-6\-5\-3\-1\-\\\end{minipage} & \Phi_{1},\Phi_{1}\Phi_{3}\Phi_{5}\\
D_8(a_2) & \begin{minipage}[t]{8em}3\-4\-2\-3\-5\-4\-2\-3\-4\-6\-5\-4\-2\-3\-7\-6\-5\-4\-2\-3\-1\-4\-5\-6\-7\-8\-\\\end{minipage} & \Phi_{2},\Phi_{2}\Phi_{6}\Phi_{10}\\
E_8(a_5) & \begin{minipage}[t]{8em}8\-6\-5\-4\-2\-3\-4\-5\-6\-7\-1\-3\-4\-5\-2\-4\-\\\end{minipage} & \Phi_{15}\\
E_8 & \begin{minipage}[t]{8em}7\-8\-6\-2\-1\-3\-4\-5\-\\\end{minipage} & \Phi_{30}\\
A_1 & \begin{minipage}[t]{8em}3\-\\\end{minipage} & \Phi_{1},\Phi_{1},\Phi_{1},\Phi_{1},\Phi_{1},\Phi_{1},\Phi_{1}\Phi_{2}\\
7A_1 & \begin{minipage}[t]{8em}7\-6\-7\-5\-6\-7\-4\-5\-6\-7\-2\-4\-5\-6\-7\-3\-4\-5\-6\-7\-2\-4\-5\-6\-3\-4\-5\-2\-4\-3\-1\-3\-4\-5\-6\-7\-2\-4\-5\-6\-3\-4\-5\-2\-4\-3\-1\-3\-4\-5\-6\-7\-2\-4\-5\-6\-3\-4\-5\-2\-4\-3\-1\-\\\end{minipage} & \Phi_{2},\Phi_{2},\Phi_{2},\Phi_{2},\Phi_{2},\Phi_{2},\Phi_{1}\Phi_{2}\\
3A_1 & \begin{minipage}[t]{8em}8\-6\-1\-\\\end{minipage} & \Phi_{1},\Phi_{1},\Phi_{1}\Phi_{2},\Phi_{1}\Phi_{2},\Phi_{1}\Phi_{2}\\
5A_1 & \begin{minipage}[t]{8em}7\-5\-4\-5\-2\-4\-5\-3\-4\-5\-2\-4\-3\-\\\end{minipage} & \Phi_{2},\Phi_{2},\Phi_{1}\Phi_{2},\Phi_{1}\Phi_{2},\Phi_{1}\Phi_{2}\\
\end{array}\]
\[
\begin{array}{lll}
\textrm{}&\textrm{}&\textrm{cont. }E_8(q)\\
\hline
A_3 & \begin{minipage}[t]{8em}4\-1\-3\-\\\end{minipage} & \Phi_{1},\Phi_{1},\Phi_{1},\Phi_{1},\Phi_{1}\Phi_{2}\Phi_{4}\\
A_3+4A_1 & \begin{minipage}[t]{8em}2\-3\-4\-2\-3\-4\-5\-4\-2\-3\-4\-5\-6\-5\-4\-2\-3\-4\-5\-6\-7\-6\-5\-4\-2\-3\-4\-5\-6\-7\-8\-\\\end{minipage} & \Phi_{2},\Phi_{2},\Phi_{2},\Phi_{2},\Phi_{1}\Phi_{2}\Phi_{4}\\
A_3+2A_1' & \begin{minipage}[t]{8em}2\-3\-4\-2\-3\-4\-5\-4\-2\-3\-4\-5\-6\-\\\end{minipage} & \left\{\begin{array}{ll}
\Phi_{1}\Phi_{2},\Phi_{1}\Phi_{2},\Phi_{1}\Phi_{2}\Phi_{4} & \textrm{ if } q \equiv 0 \mod 2\\
2,2,1/2 \Phi_{1}\Phi_{2},1/2 \Phi_{1}\Phi_{2},\Phi_{1}\Phi_{2}\Phi_{4} & \textrm{ if } q \equiv 1 \mod 2\\
\end{array}\right.
\\
 & & \textrm{uniform: }\Phi_{1},\Phi_{1},\Phi_{2},\Phi_{2},\Phi_{1}\Phi_{2}\Phi_{4}\\
D_4(a_1)+A_1 & \begin{minipage}[t]{8em}8\-5\-4\-2\-3\-4\-2\-\\\end{minipage} & \left\{\begin{array}{ll}
\Phi_{1},\Phi_{1}\Phi_{4},\Phi_{1}\Phi_{2}\Phi_{4} & \textrm{ if } q \equiv 0 \mod 2\\
2 \Phi_{1},\Phi_{1}\Phi_{4},1/2 \Phi_{1}\Phi_{2}\Phi_{4} & \textrm{ if } q \equiv 1 \mod 2\\
\end{array}\right.
\\
 & & \textrm{uniform: }\Phi_{1}\Phi_{2},\Phi_{1}\Phi_{4},\Phi_{1}\Phi_{4}\\
2A_3+A_1 & \begin{minipage}[t]{8em}2\-4\-2\-3\-4\-5\-4\-2\-3\-1\-4\-3\-5\-6\-5\-4\-2\-3\-4\-5\-6\-7\-6\-5\-4\-2\-3\-1\-4\-3\-5\-6\-7\-\\\end{minipage} & \left\{\begin{array}{ll}
\Phi_{2},\Phi_{2}\Phi_{4},\Phi_{1}\Phi_{2}\Phi_{4} & \textrm{ if } q \equiv 0 \mod 2\\
2 \Phi_{2},\Phi_{2}\Phi_{4},1/2 \Phi_{1}\Phi_{2}\Phi_{4} & \textrm{ if } q \equiv 1 \mod 2\\
\end{array}\right.
\\
 & & \textrm{uniform: }\Phi_{1}\Phi_{2},\Phi_{2}\Phi_{4},\Phi_{2}\Phi_{4}\\
D_4(a_1)+A_3 & \begin{minipage}[t]{8em}2\-3\-4\-5\-4\-2\-3\-4\-5\-7\-6\-8\-7\-6\-5\-4\-2\-3\-4\-5\-6\-\\\end{minipage} & \Phi_{4},\Phi_{4},\Phi_{1}\Phi_{2}\Phi_{4}\\
A_3+2A_1'' & \begin{minipage}[t]{8em}8\-5\-2\-4\-1\-\\\end{minipage} & \Phi_{1}\Phi_{2},\Phi_{1}\Phi_{2},\Phi_{1}\Phi_{2}\Phi_{4}\\
A_2+A_1 & \begin{minipage}[t]{8em}4\-2\-1\-\\\end{minipage} & \Phi_{1},\Phi_{1},\Phi_{1},\Phi_{1},\Phi_{1}\Phi_{2}\Phi_{3}\\
D_4+3A_1 & \begin{minipage}[t]{8em}2\-3\-4\-2\-3\-4\-5\-4\-2\-3\-4\-5\-6\-5\-4\-2\-3\-4\-5\-6\-7\-6\-5\-4\-2\-3\-1\-4\-5\-6\-7\-\\\end{minipage} & \Phi_{2},\Phi_{2},\Phi_{2},\Phi_{2},\Phi_{1}\Phi_{2}\Phi_{6}\\
3A_2+A_1 & \begin{minipage}[t]{8em}8\-5\-6\-4\-5\-6\-3\-4\-5\-6\-2\-4\-3\-1\-3\-4\-5\-6\-2\-4\-5\-3\-4\-2\-1\-\\\end{minipage} & \Phi_{3},\Phi_{3},\Phi_{1}\Phi_{2}\Phi_{3}\\
E_7(a_4) & \begin{minipage}[t]{8em}5\-2\-4\-5\-1\-3\-4\-5\-6\-2\-4\-5\-3\-4\-2\-1\-3\-4\-5\-6\-7\-\\\end{minipage} & \Phi_{6},\Phi_{6},\Phi_{1}\Phi_{2}\Phi_{6}\\
A_2+3A_1 & \begin{minipage}[t]{8em}8\-5\-2\-3\-1\-\\\end{minipage} & \Phi_{1}\Phi_{2},\Phi_{1}\Phi_{2},\Phi_{1}\Phi_{2}\Phi_{3}\\
D_4+A_1 & \begin{minipage}[t]{8em}7\-5\-4\-2\-3\-\\\end{minipage} & \Phi_{1}\Phi_{2},\Phi_{1}\Phi_{2},\Phi_{1}\Phi_{2}\Phi_{6}\\
2A_2+A_1 & \begin{minipage}[t]{8em}8\-6\-5\-3\-1\-\\\end{minipage} & \Phi_{1},\Phi_{1}\Phi_{3},\Phi_{1}\Phi_{2}\Phi_{3}\\
D_6(a_2)+A_1 & \begin{minipage}[t]{8em}2\-3\-4\-2\-3\-1\-4\-5\-4\-2\-6\-5\-4\-2\-3\-1\-4\-3\-5\-4\-6\-5\-7\-\\\end{minipage} & \Phi_{2},\Phi_{2}\Phi_{6},\Phi_{1}\Phi_{2}\Phi_{6}\\
A_5+A_2 & \begin{minipage}[t]{8em}1\-2\-3\-1\-4\-2\-3\-5\-4\-2\-3\-1\-4\-3\-5\-4\-6\-5\-4\-2\-3\-4\-5\-6\-7\-\\\end{minipage} & \left\{\begin{array}{ll}
\Phi_{3},\Phi_{1}\Phi_{2}\Phi_{3}\Phi_{6} & \textrm{ if } q \equiv 0,2 \mod 3\\
3,\Phi_{3},1/3 \Phi_{1}\Phi_{2}\Phi_{3}\Phi_{6} & \textrm{ if } q \equiv 1 \mod 3\\
\end{array}\right.
\\
 & & \textrm{uniform: }\Phi_{3},\Phi_{3},\Phi_{1}\Phi_{2}\Phi_{6}\\
E_6(a_2)+A_1 & \begin{minipage}[t]{8em}1\-2\-3\-4\-2\-6\-5\-4\-2\-3\-4\-5\-8\-\\\end{minipage} & \left\{\begin{array}{ll}
\Phi_{6},\Phi_{1}\Phi_{2}\Phi_{3}\Phi_{6} & \textrm{ if } q \equiv 0,1 \mod 3\\
3,\Phi_{6},1/3 \Phi_{1}\Phi_{2}\Phi_{3}\Phi_{6} & \textrm{ if } q \equiv 2 \mod 3\\
\end{array}\right.
\\
 & & \textrm{uniform: }\Phi_{6},\Phi_{6},\Phi_{1}\Phi_{2}\Phi_{3}\\
A_5 & \begin{minipage}[t]{8em}7\-5\-2\-6\-4\-\\\end{minipage} & \Phi_{1},\Phi_{1},\Phi_{1}\Phi_{2}\Phi_{3}\Phi_{6}\\
\end{array}\]
\[
\begin{array}{lll}
\textrm{}&\textrm{}&\textrm{cont. }E_8(q)\\
\hline
A_5+2A_1 & \begin{minipage}[t]{8em}8\-2\-3\-4\-2\-3\-4\-5\-4\-2\-3\-1\-4\-5\-6\-\\\end{minipage} & \Phi_{2},\Phi_{2},\Phi_{1}\Phi_{2}\Phi_{3}\Phi_{6}\\
D_5 & \begin{minipage}[t]{8em}2\-3\-5\-1\-4\-\\\end{minipage} & \Phi_{1},\Phi_{1},\Phi_{1}\Phi_{2}\Phi_{8}\\
D_5+2A_1 & \begin{minipage}[t]{8em}2\-3\-4\-2\-3\-4\-5\-4\-2\-3\-4\-5\-7\-6\-8\-\\\end{minipage} & \Phi_{2},\Phi_{2},\Phi_{1}\Phi_{2}\Phi_{8}\\
A_7' & \begin{minipage}[t]{8em}2\-3\-5\-6\-5\-4\-2\-3\-4\-7\-6\-5\-4\-2\-3\-1\-4\-\\\end{minipage} & \left\{\begin{array}{ll}
\Phi_{1}\Phi_{2}\Phi_{4}\Phi_{8} & \textrm{ if } q \equiv 0 \mod 2\\
2,2,1/4 \Phi_{1}\Phi_{2}\Phi_{4}\Phi_{8} & \textrm{ if } q \equiv 1 \mod 2\\
\end{array}\right.
\\
 & & \textrm{uniform: }\Phi_{4},\Phi_{8},\Phi_{1}\Phi_{2}\\
A_7'' & \begin{minipage}[t]{8em}4\-6\-8\-1\-3\-5\-7\-\\\end{minipage} & \Phi_{1}\Phi_{2}\Phi_{4}\Phi_{8}\\
A_4+A_1 & \begin{minipage}[t]{8em}5\-6\-2\-4\-1\-\\\end{minipage} & \Phi_{1},\Phi_{1},\Phi_{1}\Phi_{2}\Phi_{5}\\
D_6+A_1 & \begin{minipage}[t]{8em}6\-7\-5\-4\-5\-2\-4\-5\-1\-3\-4\-5\-2\-4\-3\-\\\end{minipage} & \Phi_{2},\Phi_{2},\Phi_{1}\Phi_{2}\Phi_{10}\\
A_3+A_2+2A_1 & \begin{minipage}[t]{8em}2\-3\-4\-2\-3\-4\-5\-4\-2\-3\-1\-4\-5\-8\-7\-\\\end{minipage} & \Phi_{2},\Phi_{2},\Phi_{1}\Phi_{2}\Phi_{3}\Phi_{4}\\
D_5(a_1) & \begin{minipage}[t]{8em}4\-2\-5\-4\-2\-3\-1\-\\\end{minipage} & \Phi_{1},\Phi_{1},\Phi_{1}\Phi_{2}\Phi_{4}\Phi_{6}\\
A_3+A_2 & \begin{minipage}[t]{8em}7\-6\-4\-2\-3\-\\\end{minipage} & \Phi_{1},\Phi_{1},\Phi_{1}\Phi_{2}\Phi_{3}\Phi_{4}\\
D_4+A_3 & \begin{minipage}[t]{8em}2\-3\-4\-2\-3\-4\-5\-7\-6\-5\-4\-2\-3\-4\-5\-6\-8\-\\\end{minipage} & \Phi_{2},\Phi_{2},\Phi_{1}\Phi_{2}\Phi_{4}\Phi_{6}\\
D_5(a_1)+A_2 & \begin{minipage}[t]{8em}1\-4\-2\-5\-4\-2\-3\-7\-8\-\\\end{minipage} & \Phi_{1}\Phi_{2}\Phi_{3}\Phi_{4}\Phi_{6}\\
D_7 & \begin{minipage}[t]{8em}6\-7\-8\-5\-4\-2\-3\-\\\end{minipage} & \Phi_{1}\Phi_{2}\Phi_{4}\Phi_{12}\\
E_6+A_1 & \begin{minipage}[t]{8em}8\-6\-4\-1\-2\-3\-5\-\\\end{minipage} & \Phi_{1}\Phi_{2}\Phi_{3}\Phi_{12}\\
E_7(a_2) & \begin{minipage}[t]{8em}4\-2\-5\-4\-3\-1\-6\-5\-7\-6\-5\-\\\end{minipage} & \Phi_{1}\Phi_{2}\Phi_{6}\Phi_{12}\\
A_6+A_1 & \begin{minipage}[t]{8em}8\-6\-7\-2\-4\-5\-1\-\\\end{minipage} & \Phi_{1}\Phi_{2}\Phi_{7}\\
E_7(a_1) & \begin{minipage}[t]{8em}5\-3\-4\-5\-6\-7\-2\-4\-1\-\\\end{minipage} & \Phi_{1}\Phi_{2}\Phi_{14}\\
E_6(a_1)+A_1 & \begin{minipage}[t]{8em}8\-5\-6\-3\-4\-2\-1\-3\-4\-\\\end{minipage} & \Phi_{1}\Phi_{2}\Phi_{9}\\
E_7 & \begin{minipage}[t]{8em}6\-7\-5\-2\-4\-3\-1\-\\\end{minipage} & \Phi_{1}\Phi_{2}\Phi_{18}\\
A_4+A_3 & \begin{minipage}[t]{8em}8\-6\-7\-4\-1\-2\-3\-\\\end{minipage} & \Phi_{1}\Phi_{2}\Phi_{4}\Phi_{5}\\
D_7(a_1) & \begin{minipage}[t]{8em}7\-8\-5\-6\-4\-2\-4\-3\-4\-\\\end{minipage} & \Phi_{1}\Phi_{2}\Phi_{4}\Phi_{10}\\
D_5+A_2 & \begin{minipage}[t]{8em}7\-8\-4\-5\-2\-3\-1\-\\\end{minipage} & \Phi_{1}\Phi_{2}\Phi_{3}\Phi_{8}\\
D_7(a_2) & \begin{minipage}[t]{8em}4\-3\-5\-4\-2\-3\-4\-5\-6\-8\-7\-\\\end{minipage} & \Phi_{1}\Phi_{2}\Phi_{6}\Phi_{8}\\
A_4+A_2+1 & \begin{minipage}[t]{8em}8\-6\-7\-5\-2\-3\-1\-\\\end{minipage} & \Phi_{1}\Phi_{2}\Phi_{3}\Phi_{5}\\
E_7(a_3) & \begin{minipage}[t]{8em}2\-4\-2\-3\-1\-6\-5\-4\-2\-3\-4\-5\-7\-\\\end{minipage} & \Phi_{1}\Phi_{2}\Phi_{6}\Phi_{10}\\
\end{array}
\]

\newpage
\subsection{Suzuki groups ${}^2B_2(q)$, $q^2 = 2^{2m+1}$}

We  need   a  factorization   $\Phi_8(q)  =  \Phi_8'(q)   \Phi_8''(q)  \in
\Q[\sqrt{2}][q]$, $\Phi_8'(q) = q^2 + \sqrt{2}  q + 1$, $\Phi_8''(q) = q^2
- \sqrt{2} q + 1$ (see~\cite[p.490]{Ca85}).

\[
\begin{array}{ll}
\textrm{class representative}&\textrm{elementary divisors}\\
\hline
-  & \Phi_{1}\Phi_{2}\\
1 & \Phi''_{8}\\
121 & \Phi'_{8}\\

\end{array}\]

\subsection{Ree groups ${}^2G_2(q)$, $q^2 = 3^{2m+1}$}

Here $\Phi_{12}'(q) =  q^2 - \sqrt{3} q  + 1$ and $\Phi_{12}''(q)  = q^2 +
\sqrt{3} q + 1$ (see~\cite[p.489]{Ca85}).

\[
\begin{array}{ll}
\textrm{class representative}&\textrm{elementary divisors}\\
\hline
- & \Phi_{1}\Phi_{2}\\
1 & \Phi'_{12}\\
121 & 2,1/2 \Phi_{4}\\
12121 & \Phi''_{12}\\

\end{array}\]

\subsection{Ree groups ${}^2F_4(q)$, $q^2 = 2^{2m+1}$}

Here  $\Phi'_{8}$   and  $\Phi''_{8}$   are  as  for   type  ${}^2B_2(q)$,
and  $\Phi_{24}'(q)  =  q^4  +\sqrt{2}  q^3  + q^2  +  \sqrt{2}  q  +  1$,
$\Phi_{24}''(q)  =   q^4  -\sqrt{2}  q^3   +  q^2   -  \sqrt{2}  q   +  1$
(see~\cite[p.490]{Ca85}).
\[
\begin{array}{ll}
\textrm{class representative}&\textrm{elementary divisors}\\
\hline
- & \Phi_{1}\Phi_{2},\Phi_{1}\Phi_{2}\\
232 & \Phi_{1}\Phi_{2}\Phi'_{8}\\
1 & \Phi_{1}\Phi_{2}\Phi_{4}\\
1213214321 & \Phi'_{24}\\
12 & \Phi''_{24}\\
2 & \Phi_{1}\Phi_{2}\Phi''_{8}\\
12132132 & \Phi_{8}\\
1232 & \Phi_{12}\\
121321324321 & \Phi_{4},\Phi_{4}\\
121321324321324321 & \Phi'_{8},\Phi'_{8}\\
121321 & \Phi''_{8},\Phi''_{8}\\
\end{array}\]

% for arXiv include .bbl file
%\bibliographystyle{alpha}
%\bibliography{gcdZX}
\newcommand{\etalchar}[1]{$^{#1}$}

\end{document}